\documentclass[12pt,reqno]{amsart}

\usepackage{hyperref}
\usepackage{amsmath}
 \usepackage[dvips]{graphicx}
 \usepackage{amsfonts}
  \usepackage{amssymb}
      \usepackage{amsthm}
 \usepackage{enumerate}
  \usepackage{latexsym}
 \usepackage{color}

 \renewcommand{\proof}{\noindent{\bf Proof. \quad}}
 \newcommand{\eproof}{\hfill \mbox{${\square}$} \vspace{3mm}}

\newcommand{\R}{\mathbb{R}}

\newcommand{\cala}{\mathcal{A}}
\def\dif{\mathrm{Diff}}

\newtheorem {teo} {\bf Theorem\,} [section]

\newtheorem {cor} [teo]{\bf Corollary}
\newtheorem {lema} [teo] {\bf Lemma}
\newtheorem {defn} [teo] {\bf Definition}
\newtheorem{rem}[teo]{\bf Remark}
\newtheorem{ex}[teo]{\bf Example}

 \makeatletter 
\def\namedlabel#1#2{\begingroup
   \def\@currentlabel{#2}%
   \label{#1}\endgroup}

\begin{document}
 
 \title[Continuity of attractors for  $\mathcal{C}^1$ perturbations]{Continuity of attractors for  $\mathcal{C}^1$ perturbations of a smooth domain}

\author[P. S. Barbosa, A. L. Pereira ]{Pricila S. Barbosa$^{\star}$ and Ant\^onio L. Pereira$^{\diamond}$ }

\thanks{$^\diamond$Partially supported by FAPESP-Brazil grant 2016/02150-8.}

\address{Ant\^onio L. Pereira \hfill\break
Instituto de Matem\'atica e Estat\'istica \\
Universidade de S\~ao Paulo - S\~ao Paulo - Brazil}
\email{alpereir@ime.usp.br}

\address{Pricila S. Barbosa \hfill\break
Universidade Tecnol\'ogica Federal do Paran\'a - Paran\'a - Brazil }
\email{pricilabarbosa@utfpr.edu.br}

\date{}


\begin{abstract}

We consider a family of semilinear parabolic problems with nonlinear boundary conditions
\[
\left\{
\begin{aligned}
u_t(x,t)&=\Delta u(x,t) -au(x,t) + f(u(x,t)) ,\,\,\ x \in \Omega_\epsilon
 \,\,\,\mbox{and}\,\,\,\,\,\,t>0\,,
\\
\displaystyle\frac{\partial u}{\partial N}(x,t)&=g(u(x,t)), \,\, x \in \partial\Omega_\epsilon \,\,\,\mbox{and}\,\,\,\,\,\,t>0\,,
\end{aligned}
\right.
\]
where $\Omega_0 \subset \R^n$ is a smooth (at least $\mathcal{C}^2$) domain, $\Omega_{\epsilon} = h_{\epsilon}(\Omega_0)$ and $h_{\epsilon}$ is a family of diffeomorphisms converging to the identity in the $\mathcal{C}^1$-norm. Assuming suitable regularity and dissipative conditions for the nonlinearites,  we show that the problem is well posed for $\epsilon>0$ sufficiently small in a suitable scale of fractional spaces,  the associated semigroup has a global attractor $\mathcal{A}_{\epsilon}$ and the family $\{\mathcal{A}_{\epsilon}\}$  is continuous at $\epsilon = 0$.

\end{abstract}

\maketitle

 \allowdisplaybreaks

\section{Introduction} \label{intro}

  Let $\Omega=  \Omega_0 \subset \R^n$ be a $\mathcal{C}^2$ domain, $a$ a positive number, $f, g: \R \to \R$ real functions, and consider  the family of semilinear parabolic problems with nonlinear Neumann boundary conditions:
 
  \begin{equation} \label{nonlinBVP} \tag{$P_{\epsilon}$}
\begin{array}{rcl}
\left\{
\begin{array}{rcl}
u_t(x,t)&=&\Delta u(x,t) -au(x,t) + f(u(x,t)) ,\,\,\ x \in \Omega_\epsilon
 \,\,\,\mbox{and}\,\,\,\,\,\,t>0\,,
\\
\displaystyle\frac{\partial u}{\partial N}(x,t)&=&g(u(x,t)), \,\, x \in \partial\Omega_\epsilon \,\,\,\mbox{and}\,\,\,\,\,\,t>0\,,
\end{array}
\right.
\end{array}
\end{equation}
where $\Omega_{\epsilon}=  \Omega_{h_\epsilon} =h_{\epsilon}(\Omega_0)$ and   $h_{\epsilon}: \Omega_0 \to \R^n $ is a family 
 of $\mathcal{C}^m, m \geq 2$  maps  satisfying suitable conditions to be specified later.
 
 One of the central questions concerning this problem is  the existence and properties of \emph{ global attractors} since, as it is well known, they  determine the  dynamics of the entire system (see, for example \cite{Hale} or \cite{Teman}). The continuity with respect to parameters present in the equation is also of interest, since it can be seen as a desirable  property of ``robustness'' in the model. In many cases, however,  the form of the equation is fixed, so the `parameter`of interest is the domain where the problem is posed.

 The existence of a global compact attractor for the problem \ref{nonlinBVP}  has been proved in 
 \cite{COPR} and \cite{OP}, under  stronger smoothness hypotheses on the domains and  growth and dissipative conditions on the nonlinearities $f$ and $g$. 

 \par The problem of existence and continuity of global attractors for semilinear parabolic problems, with respect to change of domains has also been considered  in  \cite{AC1}, for the problem with homogeneous boundary conditions 
\begin{equation*}
\begin{array}{rlr}
\left\{
\begin{array}{rcl}
u_t = \Delta u + f(x,u) \,\,\,\mbox{in}\,\,\, \,\,\,\Omega_\epsilon \\
\displaystyle\frac{\partial u}{\partial N}= 0 \,\,\,\mbox{on}\,\,\,\,\,\,\partial\Omega_\epsilon\,,
\end{array}
\right.
\end{array}
\end{equation*}
where  $\Omega_\epsilon$, $0 \leq \epsilon \leq \epsilon_0$ are bounded domains
with  Lipschitz  boundary in $\R^N$, $N \geq 2$. There it is  proved that, if the perturbations are such that the convergence of the eigenvalues and eigenfunctions of the linear part of the problem can be shown, than the upper semicontinuity of attractors follow. With the additional assumption that the equilibria are all hyperbolic, the lower semicontinuity is also obtained.

\par  The behavior of the equilibria of (\ref{nonlinBVP}) was studied in
  \cite{AB2} and \cite{AB3}. In these papers, the authors consider
 a family of smooth domains  $\Omega_\epsilon \subset \R^N$, $N \geq 2$ and  $0 \leq \epsilon \leq \epsilon_0$ whose boundary oscillates rapidly when the parameter
 $\epsilon \to 0$ and prove that the equilibria, as well as the spectra of the linearised problem around them, converge to the solution of a ``limit problem''.
 
  In \cite{PP} the authors prove the continuity of the attractors of \ref{nonlinBVP} with respect to 
 $C^2$-perturbations of a smooth domain of $\R^n$.
  
  These results do not extend immediately to the case considered here, due to the lack of smoothness of the domains considered and the fact that the perturbations do not
 converge to the inclusion in the $C^2$-norm.

 In this work, we follow the general   approach of \cite{PP}, which consists basically in ``pull-backing'' the perturbed problems to the fixed domain $ {\Omega}$
 and then considering the  family of abstract semilinear problems thus generated. We present a brief overview of this approach in the next section for convenience. Our aim here is then to prove well-posedness, 
  establish the existence  of a global attractor $\mathcal{A}_{\epsilon}$
  for sufficiently small $\epsilon \geq 0$ and prove that  the family of attractors of
 is continuous at $ \epsilon = 0$
 
\begin{figure}[!h] \label{figura}
 \centering
 \includegraphics[width=.70\columnwidth]{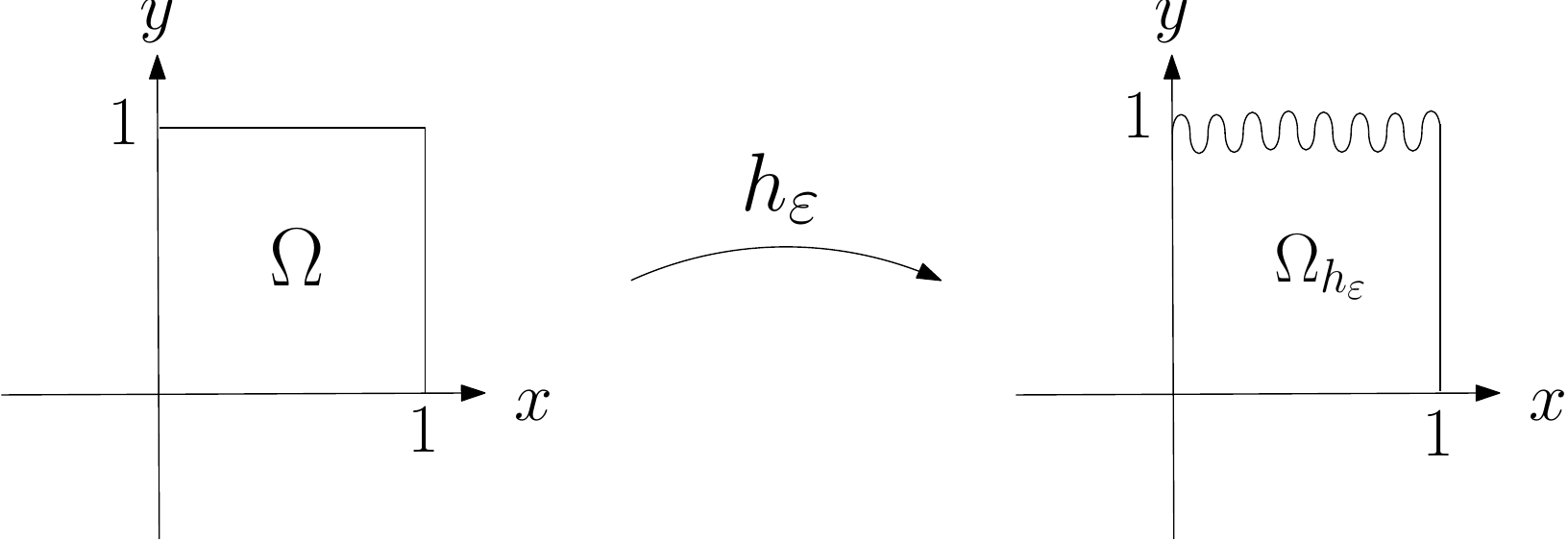}
 \caption{The perturbed region}
 \end{figure}

  These results were  obtained in our previous paper \cite{BPP} for the   family of perturbations
of the unit square  in $\R^2$  given by 

\begin{equation} \label{per}
h_{\epsilon} (x_1,x_2) =  (\,x_1\,,\,x_2 + x_2\,\epsilon\, sen(x_1/\epsilon^\alpha)\,)
\end{equation}
 with $0 <\alpha <1$  and  $\epsilon >0$  sufficiently small, 
 (see figure  (\ref{figura})).

 In the present paper, we generalize these results in two directions:  we consider the problem in arbitrary spatial dimension and, also, instead of a specific family of perturbations, we consider general  families 
   $h_{\epsilon}: \Omega_0 \to \R^n $ 
 of $\mathcal{C}^m, m \geq 2$  maps  satisfiyng the following abstract hypotheses: 
 
 \vspace{3mm}
 
  \begin{itemize}
\item  $({\bf H_1})$ \   $  \|h_{\epsilon}- i_{\Omega_0} \|_{\mathcal{C}^1{(\Omega)}} \to 0$   as   
$\epsilon \to 0.$ 
\item   $ ({\bf H_2})$ \   The Jacobian determinant $Jh_{\epsilon}$ of $h_{\epsilon}$ is differentiable, and
\\ $  \|\nabla Jh_{\epsilon} \|_{\infty} = 
 \sup \{ \, \|\nabla Jh_{\epsilon}(x) \| \ , \ x \in \Omega \}  \to 0 $ as $\epsilon \to 0$.
\end{itemize}

\vspace{3mm}

We show in section   \ref{linear_semi}  that  the family   $h_{\epsilon}$ considered in \cite{BPP} satisfies 
the  conditions $({\bf H_1})$ and $({\bf H_2})$.  Since  the domain $\Omega$ is not of class $\mathcal{C}^1$, the results obtained here do not immediately apply. However, since the perturbations occur only in a
 smooth portion of the boundary, they could easily be adapted to this case.   We also give there more general examples of families satisfying our   properties.

 The paper is organized as follows: in section \ref{prelim} we we show how the problem can be reduced to a family of problems in the initial domain   and  collect some results needed later.  In section 
 \ref{domains} we  give   some rather general examples of families  satisfying our basic assumptions.  In section 
  \ref{linear_semi} we show that the perturbed linear operators are sectorial operators in suitable spaces and  study properties of the linear semigroup generated by them.  In section \ref{abstract} we show that the problem \ref{nonlinBVP}  can be reformulated as an abstract problem in a scale of Banach spaces  which are shown to be locally well-posed in section \ref{wellposed}, under suitable growth assumptions onf
 $f$ and $g$.  In  section \ref{global_exist},  assuming a dissipative condition for the problem, we use comparison results to  prove that the solutions are globally defined and the family of associated  semigroups are uniformly bounded.
 In section \ref{attract} we prove the existence of global attractors.
 In section \ref{upper}, we show  that these attractors
 behave upper semicontinuously . Finally, in section \ref{lower}, with some additional properties on the nonlinearities and on the set of equilibria, we  show that they are also lower semicontinuous at $\epsilon=0$.

 \section{Reduction to a fixed domain} \label{prelim}

 One of the difficulties encountered in problems of perturbation of the domain
 is   that the function spaces change with the change of the region. One way to  overcome this difficulty is to  effect a ``change of variables'' in order to bring
 the problem back to a fixed region. This approach was developed by D. Henry in \cite{He1} and is the one we adopt here. 
 We describe it briefly here, for convenience of the reader. 
For a different
 approach, see \cite{AB2}, \cite{AB3} and \cite{AC1}.

 Given an open  bounded $C^m$  region $\Omega  \subset \R^{n}$, $m \geq 1$,  denote by $\textrm{Diff}^m(\Omega),
 m \geq 0$, 
the set of $\mathcal{C}^m $ embeddings (=diffeomorphisms from $\Omega$
 to  its  image).   

We define a topology in  $\textrm{Diff}^m(\Omega)$, by declaring that $ \Omega$ is in a $\epsilon$ neighborhood of $\Omega_0$, if  $ \Omega = h(\Omega_0)$, with  $ \|h\|_{\mathcal{C}^m(\Omega_0 )} < \epsilon $.   
 It has been shown   in  \cite{Mich} that this topolgy is metrizable and  we 
 denote by 
${\mathcal{M}}_{m}(\Omega)$  or simply ${\mathcal{M}}_{m}$ this (separable)
metric space.
 We say that a function $F$ defined in the
space ${\mathcal{M}}_{m}$ with values in a Banach space is $C^m$ or
analytic if $h \mapsto F(h(\Omega))$ is $C^m$ or analytic as a map
of Banach spaces  ($h$ near $i_{\Omega}$ in $C^m(\Omega,
\R^n)$). In this sense, we may express problems of perturbation
of the boundary of a boundary value problem as problems of
differential calculus in Banach spaces.

If $h: \Omega \mapsto \R^n$ is a $C^k$, $k \leq m$ embedding, we may
consider the `pull-back' of $h$

$$ h^{*}: C^k(h(\Omega)) \to C^k(\Omega) \quad (0 \leq k   \leq m)$$
defined by $h^{*}(\varphi) = \varphi \circ h$, which is an
isomorphism with inverse ${h^{-1}}^{*}$. Other function spaces can
be used instead of $C^k$, and we will actually be interested
 mainly in Sobolev spaces and fractional power spaces. 

Now, if 
 $
F_{h(\Omega)} : C^m(h(\Omega)) \to C^0(h(\Omega))$ is a (generally
nonlinear) differential operator in $\Omega_h= h(\Omega) $
   we  may consider the operator 
$h^{*} F_{h(\Omega)} {h^{*}}^{-1}$, which is a differential
operator in the fixed region $\Omega$.

%
%
%
%
%
%
%

Let  now  $h_{\epsilon}: \Omega_0 \to \R^n $ be a family of maps satisfying
 the conditions $({\bf H_1})$ and $({\bf H_2})$ and 
  $\Omega_{\epsilon} = h_{\epsilon}(\Omega)$ the corresponding family of ``perturbed domains''.

\begin{lema}\label{estimateh}
If  $\epsilon> 0$ is sufficiently small, the map   $h_\epsilon$  belongs to    $\textrm{Diff}^m(\Omega) = $  diffeomorphisms from  
 $\Omega$  { to its image}.
\end{lema}
\proof  Straightforward.

\begin{lema} \label{isosobolev}
If  $  0< s \leq m$
 and   $\epsilon > 0$ is small enough, the map 
$$
\begin{array}{llcll}
h_\epsilon^* &:& H^s(\Omega_{\epsilon}) &\to& H^s(\Omega) \\
&& u & \longmapsto & u \circ h_\epsilon
\end{array}
$$
is an isomorphism, with inverse ${h_\epsilon^*}^{-1} = (h_\epsilon^{-1})^*$.
\end{lema}

\proof See \cite{BPP}.
 \eproof
 
Using d Lemma\ref{estimateh}  we may bring  the problem \ref{nonlinBVP} back  to the fixed region $\Omega_0$.
 For this purpose, observe that 
$v(.\,,t)$  is a solution  (\ref{nonlinBVP}) in  the perturbed region
 ${\Omega_{\epsilon} = h_\epsilon(\Omega)}$, if and only if   
 $u(.\,,t)={h_\epsilon^*v(.,t)}$  satisfies
\begin{equation}\label{nonlinBVP_fix}
\begin{array}{rcl}
\left\{
\begin{array}{rcl}
u_t(x,t)&=& {h_\epsilon^*\Delta_{\Omega_{\epsilon}}h_\epsilon^{*^{-1}}} u(x,t) -au(x,t) + f(u(x,t)) ,\,\, x \in \Omega \,\,\,\mbox{and}\,\,\,t>0, \\
{h_\epsilon^*\displaystyle\frac{\partial }{\partial N_{\Omega_{\epsilon}}}h_\epsilon^{*^{-1}}}u(x,t)&=&g(u(x,t)), \,\, x \in \partial\Omega \,\,\,\mbox{and}\,\,\,t>0\,,
\end{array}
\right.
\end{array}
\end{equation}
where $ h_{\epsilon}^{*} \Delta_{\Omega_{\epsilon}} {h_{\epsilon}^{*}}^{-1}$ and
 $  h_{\epsilon}^{*} \frac{\partial}{\partial N_{\Omega_{\epsilon}}} {h_{\epsilon}^{*}}^{-1}$
are defined by
$$  
h_{\epsilon}^{*} \Delta_{\Omega_{\epsilon}} {h_{\epsilon}^{*}}^{-1}u (x) =
\Delta_{\Omega_{\epsilon}} (u\circ h_{\epsilon}^{-1})(h_{\epsilon}(x))
$$ 
and 
$$
h_{\epsilon}^{*} \frac{\partial}{\partial N_{\Omega_{\epsilon}}} {h_{\epsilon}^{*}}^{-1}=
 \frac{\partial}{\partial N_{\Omega_{\epsilon}}}(u\circ h^{-1})(hV(x))
$$
(in appropriate spaces). In particular, if
${\mathcal{A}}_{\epsilon}$ is the global attractor of (\ref{nonlinBVP}) in $
H^s(\Omega_{\epsilon})$, then $ {\tilde{\mathcal{A}}}_{\epsilon} =  \{v \circ h \mid  v
\in {\mathcal{A}}_{\epsilon}  \}$ is the global attractor of (\ref{nonlinBVP_fix})
in  $H^s(\Omega)$ and conversely. In this way we can consider
the problem of continuity of the attractors as $\epsilon \to
 0 $ in a fixed phase space.

 For later use, we   now compute an expression for the
  differential operator  $h_\epsilon^*\Delta_{\Omega_{\epsilon}}h_\epsilon^{*-1}$ in the fixed region  $\Omega$, in terms of $h_\epsilon$. 
  
  \noindent Writing
 $h_\epsilon(x)\,=\,h_\epsilon(x_1,x_2, \cdots, x_n)\,=\,((h_\epsilon)_1(x),(h_\epsilon)_2(x), \cdots, (h_\epsilon)_n(x) )\,=\,(y_1,y_2, \cdots, y_n)\,=\,y\,$, we obtain, for  $i=1, 2, \cdots, n$
\begin{equation}
\begin{split} \label{deriv}
\left(h_\epsilon^* \displaystyle\frac{\partial}{\partial y_i} h_\epsilon^{*-1}(u)\right)(x) &=  \displaystyle\frac{\partial}{\partial y_i}(u\circ h_\epsilon^{-1})(h_\epsilon(x))  \\
&= \displaystyle\sum^n_{j=1} \left[\left(\displaystyle\frac{\partial h_\epsilon}{\partial {x_j}}\right)^{-1}\right]_{j,i}(x)\frac{\partial u}{\partial x_j}(x) \\
&= \displaystyle\sum^{{n}}_{j=1} b_{ij}^{\epsilon}(x)\displaystyle\frac{\partial u}{\partial x_j}(x)\,,
\end{split}
\end{equation}

where  $b_{ij}^{\epsilon}(x)$ is the   $i,j$-entry of the
inverse transpose of the  Jacobian matrix  of  $h_\epsilon$. From now on, we omit the
$\epsilon$ from the notation for simplicity.  
  Therefore,
\begin{equation}
 \begin{split} \label{laph}
 h_\epsilon^*\Delta_{\Omega_{\epsilon}}h_\epsilon^{*-1}(u)(x)  & =   
 \sum_{i=1}^n \left(h_\epsilon^*\frac{\partial^2}{\partial y_i^2}h_\epsilon^{*-1}(u)\right)(x)   \\
 & =    \sum_{i=1}^n \left(\sum_{k=1}^n b_{i\,k} \frac{\partial}{\partial x_k}\left(\sum_{j=1}^n b_{i j}\frac{\partial u}{\partial x_j}\right) \right)(x)   \\
 & =    \sum_{k=1}^n  \frac{\partial}{\partial x_k}
 \left( \sum_{j=1}^n \sum_{i=1}^n b_{i j} b_{i k}    
 \frac{\partial u}{\partial x_j} \right){(x)}  -
 \sum_{j=1}^n \left( \sum_{i,k=1}^n
 \frac{\partial}{\partial x_k} (b_{i k} )  b_{i j} \right)
 \frac{\partial u}{\partial x_j} {(x)}\\
 & =    \sum_{k=1}^n  \frac{\partial}{\partial x_k}
 \left( \sum_{j=1}^n C_{kj}    
 \frac{\partial u}{\partial x_j} \right){(x)} -
 \sum_{j=1}^n  A_j 
 \frac{\partial u}{\partial x_j}{(x)}, 
  \end{split}
\end{equation}
where $ C_{kj} = \sum_{i=1}^n b_{i j} b_{i k} $ and
$A_j = \sum_{i,k=1}^n
 \frac{\partial}{\partial x_k} (b_{i k} )  b_{i j} $.

 We also need to compute the boundary condition $h_\epsilon^*\displaystyle\frac{\partial}{\partial N_{\Omega_{\epsilon}}}h_\epsilon^{*-1}u=0$  in the fixed region $\Omega$ in terms of  $h_\epsilon$.  Let  $N_{h_\epsilon(\Omega)}$ denote the {outward unit} normal to the boundary of  $h_\epsilon(\Omega):=\Omega_{\epsilon}$. From  (\ref{deriv}), we obtain 
\begin{equation}
\begin{split}
\left(h_\epsilon^*\displaystyle\frac{\partial}{\partial N_{\Omega_{\epsilon}}}h_\epsilon^{*-1}u\right)(x) &=  \sum_{i=1}^n
\left(h_\epsilon^*\displaystyle\frac{\partial}{\partial y_i}h_\epsilon^{*-1}u\right)(x)\left(N_{\Omega_{\epsilon}}\right)_i(h_\epsilon (x))   \\
&=  \sum_{i=1}^n \displaystyle\frac{\partial}{\partial y_i}(u \circ h_\epsilon^{-1})(h_\epsilon (x))\left(N_{\Omega_{\epsilon}}\right)_i(h_\epsilon (x))   \\
&=  \sum_{i,j=1}^n b_{ij}(x)\displaystyle\frac{\partial u}{\partial x_j}(x)\left(N_{\Omega_{\epsilon}}\right)_i(h_\epsilon(x))   \label{normalh}
\end{split}
\end{equation}
\par Since
$$ N_{\Omega_{\epsilon}}(h_\epsilon(x)) = h_\epsilon^*N_{\Omega_{\epsilon}}(x) =
\displaystyle\frac{[h_\epsilon^{-1}]_x^T N_\Omega(x)}{||\,[h_\epsilon^{-1}]_x^T N_\Omega(x)\,||} 
 $$
(see \cite{He1}), we obtain

$$ \left(N_{\Omega_{\epsilon}}(h_\epsilon(x)) \right)_i =
\frac{1}{||\,[h_\epsilon^{-1}]_x^T N_\Omega(x)\,||}
\sum_{k=1}^n b_{ik} (N_{\Omega})_k.$$

Thus, from (\ref{normalh})

\begin{equation}
\begin{split}
  \left(h_\epsilon^*\displaystyle\frac{\partial}{\partial N_{\Omega_{\epsilon}}}h_\epsilon^{*-1}u\right)(x) &=
  \frac{1}{||\,[h_\epsilon^{-1}]_x^T N_\Omega(x)\,||}
  \sum_{k=1}^n  \left(  \sum_{i,j=1}^n  b_{ik} b_{ij}(x)\displaystyle\frac{\partial u}{\partial x_j}(x) \right)   (N_{\Omega})_k \\
  &= \frac{1}{||\,[h_\epsilon^{-1}]_x^T N_\Omega(x)\,||}
  \sum_{k=1}^n  \left(  \sum_{j=1}^n   C_{kj}\displaystyle\frac{\partial u}{\partial x_j}(x) \right)   (N_{\Omega})_k 
   \end{split}
\end{equation}

 Thus, the boundary condition
$
\left(h_\epsilon^*\displaystyle\frac{\partial}{\partial N_{\Omega_{\epsilon}}}h_\epsilon^{*-1}u\right)(x) = 0\,, 
$
becomes

$$
\sum_{j,k=1}^n \left(N_\Omega(x)\right)_k(C_{kj}D_ju) = 0 \,\,\mbox{on}\,\, \partial\Omega\,.
$$

%
\noindent  so the boundary condition is exactly the
``oblique normal derivative''  with respect to the divergence part of the
 operator 
$ h_\epsilon^*\Delta_{\Omega_{\epsilon}}h_\epsilon^{*-1}.$


 \section{Basic assumptions  and examples  on domain perturbations} 
 \label{domains}
 
 We assume that  the unperturbed domain $ \Omega_0 $ is of class $\mathcal{C}^2$, and   consider rather general  examples of  families 
   $h_{\epsilon}: \Omega_0 \to \R^n $ 
 of $\mathcal{C}^2 $  maps  satisfying the  hypotheses  $({\bf H_1})$ and $({\bf H_2})$ stated in the introduction.  
 
%
%

 \begin{ex}

   The  family $h_{\epsilon}$  of perturbations 
of the unit square  in $\R^2$   considered in  \cite{BPP}   given by 

\begin{equation} 
h_{\epsilon} (x_1,x_2) =  (\,x_1\,,\,x_2 + x_2\,\epsilon\, sen(x_1/\epsilon^\alpha)\,)
\end{equation}
 with $0 <\alpha <1$  and  $\epsilon >0$  sufficiently small, 
 (see figure  (\ref{figura})).
 satisfy the  conditions $({\bf H_1})$ and $({\bf H_2})$. 
We observe that the unperturbed region \emph{is not of class $\mathcal{C}^2 $} and, therefore,  does not strictly satisfies our hypothesis. However, since the perturbation occurs only at a smooth portion of the boundary and the elliptic problem in this case is well posed (see \cite{Gri}) the problem can actually be included in the framework considered here, with only minor modifications.

    In fact,  hypothesis  $({\bf H_1})$  was  shown in  \cite{BPP}
 (Lemma 2.1). A simple computation gives $\nabla Jh_{\epsilon} = ( \epsilon^{(1-\alpha)}
  \cos(x_1/ \epsilon^{\alpha} ), 0) $, from which $({\bf H_2})$ follows easily.
  
   From  $({\bf H_1})$, it follows that  the  boundary Jacobian  $\mu_{\epsilon}= J_{\partial \Omega}   {h_{\epsilon}}_{|\, \partial \Omega} \to 1$ uniformly as $\epsilon \rightarrow 0$ . It can be checked by explicitly computation,  as done in \cite{BPP}:

 {\[ \mu_{\epsilon}  =  \left\{  
 \begin{array}{l} 
   \displaystyle\frac{\sqrt{1 + \epsilon^{\,2-2\alpha}{cos}^{\,2}  (x_1/\epsilon^{\alpha})}}{1 + \epsilon\sin(x_1/\epsilon^\alpha)} 
  \textrm{ for }  x \in I_1:= \left\{  (x_1,1) \, | \, 0
\leq x_1 \leq 1 \right\}, \nonumber \\
    \displaystyle\frac{1}{ 1+ \epsilon \sin(x_1/\epsilon^{\alpha})} \textrm{ for }  x \in I_3:= \left\{  (x_1,0) \, | \, 0
\leq x_1 \leq 1 \right\}, \nonumber \\
    1    \textrm{ for }  x \in I_2:= \left\{  (1,x_2) \, | \, 0
\leq x_2 \leq 1 \right\}  \textrm{  and }  I_4:= \left\{  (0,x_2) \, 
  | \, {0} \leq x_2 \leq 1 \right\}.  \nonumber \\
\end{array} \right.
\]
}


    \end{ex}
 Much  more general  families  satisfying the conditions $({\bf H_1})$ and  $({\bf H_2})$ are given in the examples below.
   
 \begin{ex}

    \end{ex}

 Let $ \Omega \subset \R^n$ be a $\mathcal{C}^{2} $ domain, and $X: U \subset \R^n  \to \R^n $ a  smooth (say $\mathcal{C}^1 $)
 vector field defined in an open set containing $\Omega$ and $x(t,x_0)$ the solution of
 \[ \left\{  
 \begin{array}{lcl}
    \frac{ dx}{dt} &  = &  X(x) \\
    x(0) &  =  & x_0.
 \end{array} \right.
 \]

 Then, the map $$ x: (t, \xi)  \mapsto x(t, \xi)  : (-r, r) \times \partial \Omega \to V \subset  \R^n $$ is a diffeomorphism for some $r>0$ and some open neighborhood $V$ of $\partial \Omega$. Let $ W$ be a (smaller) open neighborhood of $\partial \Omega$, that is, with  $\overline{W} \subset V$ and define $ h_{\epsilon} : W \to \R^n$ by $ h_{\epsilon}(x(t,\xi))  =  (x(t+ \eta(t) \cdot \theta_{\epsilon}(\xi)  ,\xi)) $ , where $\theta_{\epsilon}: \partial \Omega \to \R $ is a $\mathcal{C}^1$ function, with $\|\theta_{\epsilon}\|_{\mathcal{C}^1(\partial \Omega)}  \to 0 $,  as $\epsilon \to 0$, $\eta: [-r, r] \to [0,1] $ is a  $\mathcal{C}^2$ function, with $\eta(0) = 1 $ and  $\eta(t) = 0$ if $|t| \geq \frac{r}{2}$.
   Observe that $h_{\epsilon}$ is well defined and $ \{ h_{\epsilon},  \  0 \leq \epsilon \leq \epsilon_0   \}$ is a family of $ \mathcal{C}^{1}$ maps for $\epsilon_0$ sufficiently small,   with $\| h_{\epsilon}  - i_{ B_{r}(\partial \Omega) } \|_{\mathcal{C}^1{(W)}}
  \to 0 $ as $\epsilon \to 0.$
 We may extend $h_{\epsilon}$ to a diffeomorphism   of $\R^n$, 
 satisfying $({\bf H_1})$,  which we still write simply
 as $  h_{\epsilon} $  by  defining it as the identity outside $W$.

  If $ \phi: U \subset \R^{n-1} \to \R^n $  is a local  coordinate system  for
 $\partial \Omega$  in a neighborhood of   $x_0 \in \partial \Omega$, then the map
  $\Psi (t, y)= x(t, \phi(y))   :  (-r, r) \times U  \to \R^n $ is a $ \mathcal{C}^1$  coordinate system around the point
  $x_0 \in \R^n$ and
  $  \Psi^{-1} h_{\epsilon} \Psi (t,y) = (t + \eta(t) \theta_{\epsilon}(\phi(y)), y)$.
 By an easy computation, we find that the Jacobian of 
 $\Psi^{-1}  h_{\epsilon} \Psi$ is given by
 $J (\Psi^{-1} h_{\epsilon} \Psi(t,y)) = 1 + \eta'(t)  \theta_{\epsilon}(\phi(y)) $ and, therefore
   $J  h_{\epsilon}  (x)=  \left[ 1 + \eta'(t(x))  \theta_{\epsilon}(\phi(\pi(x))) 
 \right]\cdot  J \Psi \left( \Psi^{-1}( h_{\epsilon} (x)) \right)\cdot    J\Psi^{-1}  (x)  $
  for $x \in W $.
   Since
 $\| h_{\epsilon}  - Id_{ R^n } \|_{\mathcal{C}^1}\to 0$, the condition $({\bf H_2})$ follows.

 We can also  compute  $J_{\partial \Omega}  {h_{\epsilon}}_{|\, \partial \Omega} $, the Jacobian of 
 $  h_{\epsilon} $ restricted to  ${\partial \Omega}$.  We drop the subscript
 $\partial \Omega$  to simplify the notation. 
Note that   the  coordinate system $\Psi$ above takes $ \{0 \} \times U$
 into a neighborhood of  $x_0 \in \partial \Omega,$ and
 $ \Psi^{-1}  {h_{\epsilon}}_{|\, \partial \Omega}  \Psi (0, y)= ( \theta_{\epsilon}(\phi(y)), y)$.
 
 A straightforward computation then gives
 $J (\Psi^{-1} {h_{\epsilon}}_{|\, \partial \Omega} \Psi(0,y)) =
  \sqrt{ 1 + \|\nabla \theta_{\epsilon} (\phi(y))\|^2} $ and, therefore
   $J  {h_{\epsilon}}_{|\, \partial \Omega}  (\phi(y))= 
 \left[ \sqrt{ 1 + \|\nabla \theta_{\epsilon} (\phi(y))\|^2}
 \right] J \Psi \left( \Psi^{-1}_{\epsilon}( h_{\epsilon} (\phi(y))) \right)\cdot    J\Psi^{-1}  
(\Psi  (0,y))  $
  for $y \in
 U $, where  $ \Psi_{0}$ and $ \Psi_{\epsilon}$ denote the restriction
 of $\Psi$ to $ \{ (0, y)  \, |y \in U \}$  and 
 $ \{ ( \theta_{\epsilon}(\phi(y)) , y)  \, |y \in U \}$, respectively.  Since
 $\| h_{\epsilon}  - Id_{ R^n } \|_{\mathcal{C}^1}$  and  $\|\theta_{\epsilon}(\xi)\|_{\mathcal{C}^1{(\partial\Omega)} }  \to 0 $,  it follows that
 $J  {h_{\epsilon}}_{|\, \partial \Omega}  (\phi(y))   \to  1 $ as $\epsilon \to 0$, uniformly in $\partial \Omega$. \eproof
  
 \begin{ex}
 \end{ex}
 We can choose the vector field $X$ in the previous example as an extension of $N : \partial \Omega  \to \R^n$  the unit outward normal to $ \partial \Omega $, $t(x) = \pm \textrm{dist}(x, \partial \Omega), \quad
 ("+" \textrm{ outside},   "-" \textrm{ inside})$,
 $\phi(x) = \textrm{ the point of } \partial \Omega \textrm{ nearest
   to } x$ and
 $B_r(\partial \Omega) = \{x \in \R^n \, | \,  \textrm{dist}(x, \partial \Omega) < r \}.$

 Then, the map $ \rho: (t, \xi) \mapsto \xi + t N(\xi) : (-r, r) \times \partial \Omega \to
 B_r(\partial \Omega)$  is a diffeomorphism, for some $r> 0$, with inverse
 $x \mapsto (t(x), \pi(x))$
 (see \cite{He1}).

 Define $ h_{\epsilon} : B_r(\partial \Omega) \to \R^n$ by
  $ h_{\epsilon}(\rho(t,\xi))  = \xi + t N(\xi) +
  \eta(t) \theta_{\epsilon}(\xi) N(\xi)  =
  \rho(t,\xi) + \eta(t)\theta_{\epsilon}(\xi) N(\xi) $, where
  $\theta_{\epsilon}: \partial \Omega \to \R $ is a $\mathcal{C}^1$ function, with
  $\|\theta_{\epsilon} \|_{\mathcal{C}^1{(\partial\Omega)}}  \to 0 $  as $\epsilon \to 0$,
 $\eta: [-r, r] \to [0,1] $ is a  $\mathcal{C}^2$ function, with 
$\eta(0) = 1 $ and  $\eta(t) = 0$ if $|t| \geq \frac{r}{2}$.
  Then,  $ \{ h_{\epsilon},  \  0 \leq \epsilon \leq \epsilon_0   \}$ is a family of $ \mathcal{C}^{1}$ maps for $\epsilon_0$ sufficiently small,   with
  $\| h_{\epsilon}  - i_{ B_{r}(\partial \Omega) } \|_{\mathcal{C}^1}
  \to 0$ as $\epsilon \to 0.$
 We may extend $h_{\epsilon}$ to a diffeomorphism   of $\R^n$, 
 satisfying $({\bf H_1})$,  which we still write simply
 as $  h_{\epsilon} $  by  defining it as the identity outside $B_r(\partial
 \Omega)$.

  If $ \phi: U \subset \R^{n-1} \to \R^n $ is a local  coordinate system  for
 $\partial \Omega$  in a neighborhood of   $x_0 \in \partial \Omega$, then the map
  $\Psi (t, y)=  \phi(y) + t N( \phi(y)) = \rho(t, \phi(y)) :  (-r, r) \times U  \to \R^n $ is a $ \mathcal{C}^1$  coordinate system around the point
  $x_0 \in \R^n$ and
  $  \Psi^{-1} h_{\epsilon} \Psi (t,y) = (t + \eta(t) \theta_{\epsilon}(\phi(y)), y)$.
   The condition  $({\bf H_2})$ can now be checked as in the previous example.

 \begin{rem}
We may choose the function $\theta_{\epsilon}$ with ``oscillatory behavior'', so the example above essentialy includes the case considered in \cite{BPP}, since the perturbation there is nonzero  only in a  smooth portion of the boundary.
   \end{rem}

 \section{The linear semigroup} \label{linear_semi}

 In this section we  consider the linear semigroups generated by the family of  differential operators
$   - h_\epsilon^*\Delta_{\Omega_{\epsilon}}h_\epsilon^{*-1} {+\,aI}$, appearing in 
 (\ref{nonlinBVP_fix}).

 \subsection{Strong form in $L^p$ spaces}
 Consider the operator in $L^p(\Omega),  \ p \geq 2 $, given by 
 
  \begin{equation}\label{operator_Lp}
A_{\epsilon}:= \left(\,- h_\epsilon^*\Delta_{\Omega_{\epsilon}}h_\epsilon^{*-1} +aI\,\right)
\end{equation}
 with domain
\begin{equation}\label{domain_Lp}
D\left(A_{\epsilon}\right) = \left\{u \in W^{2,p} (\Omega) \,\bigg|\, 
  h_{\epsilon}^*\displaystyle\frac{\partial }{\partial N_{\Omega_{\epsilon}}}h_{\epsilon}^{*^{-1}}u = 0,\, \textrm{ on }  \,  \partial \Omega \right\}.
\end{equation}
 
(We will  denote simply by $A$  the unperturbed operator  
$\left(\,{-}\Delta_{\Omega}  +aI\,\right)$).

\begin{teo}\label{sectorial_Lp}
If  $\epsilon> 0$ is  sufficiently small and  $h_\epsilon \in \dif^1(\Omega)$,
then the operator $A_{\epsilon}= \left(\,- h_\epsilon^*\Delta_{\Omega_{\epsilon}}h_\epsilon^{*-1} +aI\,\right)$ defined by  (\ref{operator_Lp}) and (\ref{domain_Lp})
 is sectorial.
\end{teo}

\proof   
Consider the operator  $- \Delta_{\Omega_{\epsilon}} $ defined in $ L^p( h_{\epsilon}(\Omega))$, with
domain 

 \[ D( - \Delta_{\Omega_{\epsilon}})  =  \left\{u \in W^{2,p} (\Omega_{\epsilon}) \,\bigg|\, 
  \frac{\partial }{\partial N_{\Omega_{\epsilon}}}u = 0,\, \textrm{ on }  \,  \partial \Omega_{\epsilon} \right\},\]
where $\Omega_{\epsilon} =  h_{\epsilon}(\Omega)$. 
It is well known that $- \Delta_{\Omega_{\epsilon}}$ is sectorial, with the spectra contained in the interval
$ ]   0,  \infty) \subset \R $.  
 
 If  $\lambda \in \mathbb{C} $ and $f \in L^2(\Omega)$, we have
 \begin{align}
  &\left( \, h_\epsilon^*\Delta_{\Omega_{\epsilon}}h_\epsilon^{*-1} + \lambda I \right) u (x)  = f (x)   \\ 
    \nonumber \\
   \Leftrightarrow   &  \left( \Delta_{\Omega_{\epsilon}} + \lambda I \right)  u \circ h^{-1}_{{\epsilon}} (h_{{\epsilon}} (x))  =
    f \circ h^{-1}_{{\epsilon}} (h_{{\epsilon}} (x))  \nonumber  \\ 
   \Leftrightarrow   &  \left( \Delta_{\Omega_{\epsilon}} + \lambda I \right)  v  (y )  =
    g (y) \nonumber,
 \end{align}
 Since $ u \mapsto h^*_{{\epsilon}}  u := u \circ h_{{\epsilon}}  $ is an isomorphism from $ L^2(\Omega_{\epsilon})$ to 
 $ L^2(\Omega )$ with inverse $ ({h^{-1}_{{\epsilon}} })^{*} $, it follows that the 
first equation is uniquely solvable in  $ L^2(\Omega )$ if and only if the last equation is 
uniquely solvable in $ L^2(\Omega_{\epsilon})$ .

Suppose $\lambda $ belongs to  $ \rho( -\Delta_{\Omega_{\epsilon}} ) $, the resolvent  set of $-\Delta_{\Omega_{\epsilon}}$. Then, we have.

\begin{align*}
 \| u \|_{L^p(\Omega)}^p  & =  \int_{\Omega}  |u (x)|^p  \, d \, x \\
  & =  \int_{\Omega}  |v \circ h_{{\epsilon}}   (x)|^p  \, d \, x  \\  
  & =  \int_{\Omega_{\epsilon} }  |v (y) |^p  |J h^{-1}_{{\epsilon}}  (y) | \, d \, y  \\  
  & \leq \|Jh^{-1}_{{\epsilon}}  \|_{\infty}  \| v \|_p^p   \\
  & \leq  \|Jh^{-1}_{{\epsilon}}  \|_{\infty}   \cdot
  \| \left( \Delta_{\Omega_{\epsilon}} + \lambda I \right)^{-1} \|_{\mathcal{L}( L^p (\Omega_{\epsilon}))} \cdot 
  \| g  \|_{{L^p(\Omega_ {\epsilon})} }^p   \\
\end{align*}

 On the other hand
 
\begin{align*}
 \| g  \|_{{L^p(\Omega_ {\epsilon})} }^p  & =   \int_{\Omega_{{\epsilon}} }  |g(x)|^p  \, d \, y \\
  & =  \int_{\Omega_{{\epsilon}} }  |f \circ h^{-1}_{{\epsilon}}   (y)|^p  \, d \, y  \\  
  & =  \int_{{\Omega} } |f (x) |^p  |Jh_{{\epsilon}}  (x) | \, d \, x  \\  
  & \leq \|Jh_{{\epsilon}}  \|_{\infty}  \| f  \|_{{L^p(\Omega )} }^p     \\
\end{align*}

 It follows that 
 \[
  \| u \|_{L^p(\Omega)}^p  \leq     \|Jh_{{\epsilon}}  \|_{\infty}    \cdot  \|Jh^{-1}_{{\epsilon}}  \|_{\infty} \cdot \| \left( \Delta_{\Omega_{\epsilon}} + \lambda I \right)^{-1} \|_{\mathcal{L}( L^p (\Omega_{\epsilon}))} \cdot  \| f  \|_{{L^p(\Omega )} }^p  
 \]

 Therefore,   $ \lambda  \in  \rho(-  h_\epsilon^* \Delta_{\Omega_{\epsilon}}h_\epsilon^{*-1}  ) $ and 
 
 \begin{equation} \label{resolv_Lp} 
 \|  \left( h_\epsilon^*  \Delta_{\Omega_{\epsilon}}  h_\epsilon^{*-1}   + \lambda I \right)^{-1} \|_{\mathcal{L}( L^p (\Omega))} \leq 
  \|Jh_{{\epsilon}}  \|_{\infty}    \cdot  \|Jh^{-1}_{{\epsilon}}  \|_{\infty} \cdot \| 
   \left(   \Delta_{\Omega_{\epsilon}}    + \lambda I \right)^{-1} \|_{\mathcal{L}( L^p (\Omega_{\epsilon}))}.
\end{equation}

  Reciprocally, one can prove similarly that   $ \lambda  \in  \rho(-  h_\epsilon^* \Delta_{\Omega_{\epsilon}}h_\epsilon^{*-1}  )
   \Rightarrow  \lambda  \in  \rho(-  \Delta_{\Omega_{\epsilon}} ).  $ 

Finally, if $ B_{\epsilon} = -  \Delta_{\Omega_{\epsilon}} + a I $ is sectorial   with  $ \|(\lambda -
B_{\epsilon}^{-1}\| \leq  \displaystyle\frac{M}{|\lambda -a' |} $ for all $\lambda$ in the
sector $ S_{ a',\phi_0} = \{ \lambda  \ | \ \phi_0 \leq
|arg(\lambda-a')| \leq \pi, \lambda \neq a' \} $, for some $a'\in
\mathbb{R}$ and $0\leq \phi_0<\pi/2$, it follows from  
   \eqref{resolv_Lp} that $ A_{\epsilon} = a-  h_\epsilon^* \Delta_{\Omega_{\epsilon}}h_\epsilon^{*-1} $
   satisfies $ \|(\lambda -
A)^{-1}\| \leq  \displaystyle\frac{M'}{|\lambda -a' |} $ for all $\lambda$ in the
 sectoriality of $A_{\epsilon}$ follows from the sectoriality of $B_{\epsilon}$.
 
\eproof

 \begin{rem}
 From  \ref{sectorial_Lp} and results in  \cite{He2}, it follows that $A_\epsilon$ generates a linear analytic semigroup in $L^p(\Omega)$, for each $\epsilon \leq 0$.                                         
\end{rem}

\subsection{Weak form in $L^p$ spaces}

 One would like to prove that the operators $A_{\epsilon}$ defined by \eqref{operator_Lp} and \eqref{domain_Lp} become close to the operator $A$ as $\epsilon \to 0$  in a certain sense.  This is possible when the perturbation diffeomorphisms $h_{\epsilon}$ converge to the identity in the
  $\mathcal{C}^2$-norm   (see, for example \cite{OPP} and \cite{PP}).
 To obtain similar results here, we need to consider the problem in weaker topologies, that is, we need to extend those operators. To this end, we   now want to consider the operator 
  $A_{\epsilon} =
\left(-h_\epsilon^*\Delta_{\Omega_{\epsilon}}h_\epsilon^{*-1} + aI\right)$  as  an operator $\widetilde{A}_{\epsilon}$
 in   $(W^{1,q}(\Omega))'$ with  $D(\widetilde{A}_{\epsilon})=W^{1,p} (\Omega)$, where
 $q$ is the conjugate exponent of $p$, that is $\frac{1}{p} + \frac{1}{q}=1$. 

\par If   $u \in D(A_{\epsilon})=\left\{\,u \in W^{2,p}(\Omega)\,\,\bigg|\,\,h_\epsilon^*\displaystyle\frac{\partial}{\partial N_{\Omega_{\epsilon}}}h_\epsilon^{*-1}u=0\,\right\}$, $\psi \in W^{1,q}(\Omega)$, and
 $v = u \circ h_\epsilon^{-1}$, we obtain, integrating by parts

\begin{align} \label{weak_form_Lp}
\left\langle A_{\epsilon}u\,,\,\psi\right\rangle_{{-1,1}} &= - \displaystyle\int_\Omega (h_\epsilon^*\Delta_{\Omega_{\epsilon}}h_\epsilon^{*-1}u)(x)\,\psi(x)\,dx+ a\displaystyle\int_\Omega u(x)\psi(x)\,dx  \nonumber \\
&= - \displaystyle\int_\Omega \Delta_{\Omega_{\epsilon}} (u \circ h_\epsilon^{-1})(h_\epsilon(x))\,\psi(x)\,dx + a\displaystyle\int_\Omega u(x)\psi(x)\,dx  \nonumber \\
&= - \displaystyle\int_{\Omega_{\epsilon}} \Delta_{\Omega_{\epsilon}} v(y)\psi(h_\epsilon^{-1}(y)) \displaystyle\frac{1}{|Jh_\epsilon(h_\epsilon^{-1}(y))|}dy + a\displaystyle\int_{\Omega_{\epsilon}} u(h_\epsilon^{-1}(y))\psi(h_\epsilon^{-1}(y))\displaystyle\frac{1}{|Jh_\epsilon(h_\epsilon^{-1}(y))|}dy \nonumber \\
&=  -\displaystyle\int_{\partial \Omega_{\epsilon}} \displaystyle\frac{\partial v}{\partial N_{\Omega_{\epsilon}}}(y)\, \psi(h_\epsilon^{-1}(y)) \displaystyle\frac{1}{|\,Jh_\epsilon(h_\epsilon^{-1}(y))\,|}\,d\sigma (y)  \nonumber \\
& + \displaystyle\int_{\Omega_{\epsilon}} \nabla_{\Omega_{\epsilon}} v(y)\cdot \nabla_{\Omega_{\epsilon}} \left[ \psi(h_\epsilon^{-1}(y)) \displaystyle\frac{1}{|\,Jh_\epsilon(h_\epsilon^{-1}(y))\,|} \right] \,dy  \nonumber \\
& \,+ a\displaystyle\int_{\Omega_{\epsilon}} u(h_\epsilon^{-1}(y))\psi(h_\epsilon^{-1}(y))\,\displaystyle\frac{1}{|\,Jh_\epsilon(h_\epsilon^{-1}(y))\,|} \,dy\, \nonumber \\
&=  
\displaystyle\int_{\Omega_{\epsilon}} \nabla_{\Omega_{\epsilon}} v(y)\cdot \nabla_{\Omega_{\epsilon}}
 \left[ \psi(h_\epsilon^{-1}(y)) \displaystyle\frac{1}{|\,Jh_\epsilon(h_\epsilon^{-1}(y))\,|} \right] \,dy
   \nonumber\\
& \, + a\displaystyle\int_{\Omega_{\epsilon}} u(h_\epsilon^{-1}(y))\psi(h_\epsilon^{-1}(y))\,\displaystyle\frac{1}{|\,Jh_\epsilon(h_\epsilon^{-1}(y))\,|}\,dy\nonumber\\
 &=  
 \displaystyle\int_{\Omega} \nabla_{\Omega_{\epsilon}}v(h_\epsilon(x))\cdot \nabla_{\Omega_{\epsilon}}
 \left[ \psi \circ h_\epsilon^{-1}\displaystyle\frac{1}{|Jh_\epsilon \circ h_\epsilon^{-1}|}(h_\epsilon(x))
  \right]
   \,  |Jh_\epsilon(x)|dx \nonumber
  + a\displaystyle\int_{\Omega} u(x)\psi(x)dx\nonumber \\
&=  
\displaystyle\int_{\Omega} (h_\epsilon^*\nabla_{\Omega_{\epsilon}}h_\epsilon^{*-1}u)(x) \cdot\left[ h_\epsilon^* \nabla_{\Omega_{\epsilon}} h_\epsilon^{*-1} \displaystyle\frac{\psi}{Jh_\epsilon}\right]
(x)\,|\,Jh_\epsilon(x)\,|\,dx\,
+ a\displaystyle\int_{\Omega} u(x)\psi(x)\,dx \nonumber  \\
&=  
\displaystyle\int_{\Omega} (h_\epsilon^*\nabla_{\Omega_{\epsilon}}h_\epsilon^{*-1}u)(x) \cdot  h_\epsilon^* \nabla_{\Omega_{\epsilon}} h_\epsilon^{*-1} \displaystyle \psi 
(x) \,dx\,  + a\displaystyle\int_{\Omega} u(x)\psi(x)\,dx  \nonumber \\
 & + 
\int_{\Omega} (h_\epsilon^*\nabla_{\Omega_{\epsilon}}h_\epsilon^{*-1}u)(x) \cdot ( h_\epsilon^* \nabla_{\Omega_{\epsilon}} h_\epsilon^{*-1} \displaystyle Jh_\epsilon)(x) \, \cdot 
\frac{1}{Jh_\epsilon }  \cdot \psi 
(x) \,dx.
\end{align}

 Since  (\ref{weak_form_Lp}) is well defined for  $u \in W^{1,p}(\Omega)$, we may define   an extension
 $\widetilde{A}_{\epsilon}$ of $A_{\epsilon}$, with  domain  $W^{1,p}(\Omega)$ and values in $(W^{1,q}(\Omega))'$, by 
  \begin{align} \label{operator_weak_Lp}
  \left\langle \widetilde{A}_{\epsilon}u\,,\,\psi\right\rangle_{{-1,1}}  :=   & 
\displaystyle\int_{\Omega} (h_\epsilon^*\nabla_{\Omega_{\epsilon}}h_\epsilon^{*-1}u)(x) \cdot  h_\epsilon^* \nabla_{\Omega_{\epsilon}} h_\epsilon^{*-1} \displaystyle \psi 
(x) \,dx\,  + a\displaystyle\int_{\Omega} u(x)\psi(x)\,dx  \nonumber  \\
  + &
\int_{\Omega} (h_\epsilon^*\nabla_{\Omega_{\epsilon}}h_\epsilon^{*-1}u)(x) \cdot ( h_\epsilon^* \nabla_{\Omega_{\epsilon}} h_\epsilon^{*-1} \displaystyle Jh_\epsilon)(x) \, \cdot 
\frac{1}{Jh_\epsilon }  \cdot \psi 
(x) \,dx,
\end{align}
for any $\Psi \in (W^{1,q}(\Omega))$.

  \begin{rem}
  If $u$ is regular enough, then  $\widetilde{A} u  = A u$   implies that $u$ must satisfy the 
  boundary condition  $ h_{\epsilon}^*\displaystyle\frac{\partial }{\partial N_{\Omega_{\epsilon}}}h_{\epsilon}^{*^{-1}}u = 0 $, on $\partial \Omega$ but, since this is not well defined in 
  $ (W^{1,q}(\Omega)) $,  the domain of $\widetilde{A}$ does not incorporate this boundary condition.
  \end{rem}
   
For simplicity, we still denote this extension by $A_{\epsilon},$ whenever there is no danger of confusion. Also,  from now on, we drop the absolute value in $ |\, Jh_{\epsilon} (x) \, |$, since    the Jacobian of
  $h_\epsilon$ is positive for sufficiently small $\epsilon$. 

 We now prove the following basic inequality.


 \begin{teo}\label{fund_inequality_Lp}    $D\big(A_{\epsilon}\big) \supset D\big(A\big)$ for any $\epsilon \geq 0$ and  there exists a positive function    $\tau(\epsilon)$    such that
$$
\big|\big|\,\big(A_{\epsilon}  - A  \big)u\,\big|\big|_{W^{1,q}(\Omega)'} \leq  {\tau}(\epsilon)\big|\big|\, A  \,u\,\big|\big|_{W^{1,q}(\Omega)'} \,,
$$ 
for all  $u \in D\big( A  \big)$, with  $\displaystyle\lim_{\epsilon \to 0{^+}}  {\tau}(\epsilon)=0.$
 \end{teo}
 \proof
 The assertion about the domain is immediate.
 The  inequality is 
  equivalent to 
$$
\big|\,\big\langle\,\big(A_{\epsilon}  -A \big)u\,,\,\psi\,\big\rangle_{-1,1}\big| 
\leq \tau(\epsilon)||\,A u\,||_{(W^{1,q}(\Omega))' }||\,\psi\,||_{W^{1,q}(\Omega)}\,, 
$$
for all  $u\in W^{1,p}(\Omega) $, $\psi \in  W^{1,q}(\Omega)$, with  $\displaystyle\lim_{\epsilon \to 0{^+}} \tau(\epsilon)=0.$ 
 We have, for $\epsilon >0$.
 \begin{align}
     \left| \big\langle \,\big(A_{\epsilon}   -A   \,\big)u\,,\,\psi\,\big\rangle_{ {-1,1}}\right| &  =
    \left |\displaystyle\int_{\Omega} \big(\,h_\epsilon^*\nabla_{\Omega_{\epsilon}}h_\epsilon^{*-1}u\,\big)(x) \cdot \left[ \big(\,h_\epsilon^*\nabla_{\Omega_{\epsilon}}h_\epsilon^{*-1} \psi \,\big)(x) -
     \big( \nabla_{\Omega}  \psi \,\big)(x) \right] \, d\,x  \right| \nonumber \\
      &  +    \left |\displaystyle\int_{\Omega} \big(\,h_\epsilon^*\nabla_{\Omega_{\epsilon}}h_\epsilon^{*-1}u -  \nabla_{\Omega} u \,\big)(x) \cdot 
    \big( \nabla_{\Omega}  \psi \,\big)(x) \, d\,x  \right| \nonumber \\
  &  +   
 \left| \int_{\Omega} (h_\epsilon^*\nabla_{\Omega_{\epsilon}}h_\epsilon^{*-1}u)(x) \cdot ( h_\epsilon^* \nabla_{\Omega_{\epsilon}} h_\epsilon^{*-1} \displaystyle Jh_\epsilon)(x) \, \cdot 
\frac{1}{Jh_\epsilon }  \cdot \psi 
(x) \,dx \right|
 \end{align}

Now, writing $ |v |_p = \left( \sum_{i=1}^n |v_i|^p \right)^{\frac{1}{p}}$, $1\leq p < \infty$, \  $ |v|_{\infty} = \sup \left(  |v_i|, \ i=1,2, \cdots, n  \right)$ for the $p$-norm of the vector 
 $v=(v_1,v_2, \cdots, v_n) \in \R^n $, we observe that
\begin{align*}
\left| \,h_\epsilon^*\nabla_{\Omega_{\epsilon}}h_\epsilon^{*-1}u\,(x)  \right|_p & =
 \left( \sum_i \left| \,h_\epsilon^*   \frac{\partial }{\partial y_i} h_\epsilon^{*-1}u\,(x)  
   \right|^p \right)^{\frac{1}{p}} 
   = \left( \sum_i   \left(  \sum_j \left|  b_{i,j}^{\epsilon} (x) \frac{\partial u}{\partial x_j} (x) \right| \right)^p \right)^{\frac{1}{p}}  \\
  &\leq
 \left[ \sum_i    \left( \sum_j | ( b_{i,j}^{\epsilon})|^q (x) \right)^{\frac{p}{q}} 
 \left( \sum_j  \left( \left|\frac{\partial u}{\partial x_j} \right|  \right)^p (x)\right)   \right]^{\frac{1}{p}}  \\ 
 & \leq
 \left[  \sum_i   \left( \sum_j  \left|( b_{i,j}^{\epsilon})\right|^q  (x)\right)^{p-1} \right]^{\frac{1}{p}} |\nabla u (x) |_p   \\
 & \leq
    \| ( b^{\epsilon})\|_{\infty}  \left[  \sum_i  n^{p-1} \right]^{\frac{1}{p}} |\nabla u (x) |_p   \\
  & \leq
    n \|  b^{\epsilon}\|_{\infty}  |\nabla u (x) |_p   \\
       & \leq  B(\epsilon)  |\nabla u (x) |_p \\
 & \\
  \left| \,h_\epsilon^*\nabla_{\Omega_{\epsilon}}h_\epsilon^{*-1}u\,(x) - \nabla_{\Omega} u (x)   \right|_p & = 
  \left( \sum_i \left| \,h_\epsilon^*   \frac{\partial }{\partial y_i} h_\epsilon^{*-1}u\,(x)  - \frac{\partial u}{\partial x_i}(x) 
   \right|^p \right)^{\frac{1}{p}} \\ 
 & =   \left( \sum_i   \left(  \sum_j \left|  \left(b_{i,j}^{\epsilon} (x) - \delta_{i,j}\right) \frac{\partial u}{\partial x_j} (x) \right| \right)^p \right)^{\frac{1}{p}}  \\
 &\leq 
 \left[ \sum_i    \left( \sum_j | ( b_{i,j}^{\epsilon} - \delta_{i,j} )|^q (x) \right)^{\frac{p}{q}} 
 \left( \sum_j  \left( \left|\frac{\partial u}{\partial x_j} \right|  \right)^p (x)\right)   \right]^{\frac{1}{p}}  \\ 
 & \leq
 \left[  \sum_i   \left( \sum_j  \left|( b_{i,j}^{\epsilon}) - \delta_{i,j} \right|^q  (x)\right)^{p-1} \right]^{\frac{1}{p}} |\nabla u (x) |_p   \\
 & \leq
    \| ( b^{\epsilon}-\delta)\|_{\infty}  \left[  \sum_i  n^{p-1} \right]^{\frac{1}{p}} |\nabla u (x) |_p   \\
  & \leq
    n \|  b^{\epsilon} -\delta \|_{\infty}  |\nabla u (x) |_p   \\
                                                                                 & \leq  \eta(\epsilon)  |\nabla u (x) |_p \\
   & \\
 \frac{1}{Jh_\epsilon(x)} \left| \,h_\epsilon^*\nabla_{\Omega_{\epsilon}}h_\epsilon^{*-1}Jh_{\epsilon}\,(x)  \right|_{\infty} & = \frac{1}{Jh_\epsilon(x)}
 \ \sup_i  \left\{ \left| \,h_\epsilon^*   \frac{\partial }{\partial y_i} h_\epsilon^{*-1}Jh_{\epsilon}\,(x)  \right|
 \right\} \\
  &  =  \frac{1}{Jh_\epsilon(x)} \sup_i   \left\{   \sum_j \left|  b_{i,j}^{\epsilon} (x) \frac{\partial Jh_{\epsilon}}{\partial x_j} (x) \right| \right\}   \\
    & =  \frac{1}{Jh_\epsilon(x)}  \ \|  b_{\infty}^{\epsilon} \| \sum_j \left| \frac{\partial Jh_{\epsilon}}{\partial x_j} (x) \right|    \\
     & \leq  \frac{1}{Jh_\epsilon(x)} \ \|  b_{\infty}^{\epsilon} \|  |\nabla Jh_{\epsilon} (x)|_1     \leq  \ n \|  b_{\infty}^{\epsilon}  \| |\nabla Jh_{\epsilon} (x)|_{\infty} \\  
                                                                                 & \leq \frac{1}{Jh_\epsilon(x)}   B(\epsilon)  |\nabla Jh_{\epsilon} (x) |_{\infty} 
                                                                                   \leq \frac{1}{Jh_\epsilon(x)}
                                                                                   B(\epsilon)  \|\nabla Jh_{\epsilon}  \|_{\infty}  \\
& \leq \mu(\epsilon) \\
                                                                                 & \\
   \frac{1}{Jh_\epsilon(x)} \left| \,h_\epsilon^*\nabla_{\Omega_{\epsilon}}h_\epsilon^{*-1}Jh_{\epsilon}\,(x)  \psi (x)  \right|_{q} & = \frac{1}{Jh_\epsilon(x)}
 \ \left( \sum_i   \left| \,h_\epsilon^*   \frac{\partial }{\partial y_i} h_\epsilon^{*-1}Jh_{\epsilon}\,(x) \cdot \psi(x)  \right|^q \right)^{\frac{1}{q}}
  \\
                                                                                 & \leq
                                                                                    \frac{1}{Jh_\epsilon(x)}
                                                                                     \left| h_\epsilon^*   \nabla  h_\epsilon^{*-1}Jh_{\epsilon}\,(x) \right|_{\infty}   \ \left(
                                                                                   \sum_i   \left|  \psi(x)  \right|^q \right)^{\frac{1}{q}}
  \\
  & \leq n  \mu(\epsilon)   \psi(x) ,
   \end{align*} 
  where  $\|  b^{\epsilon} \|_{\infty};= \sup \{  | b_{i,j}^{\epsilon} |  (x),   \  \, 1 \leq i,j \leq n, \ x \in \Omega   \}, $
  $ \|  b^{\epsilon} -\delta \|_{\infty};= \sup \{  | b_{i,j}^{\epsilon} - \delta_{i,j}  |  (x),   \  \, 1 \leq i,j \leq n, \ x 
  \in \Omega   \},  $ \     $B(\epsilon) \to n $ and  $\eta(\epsilon)$,
$\mu(\epsilon)  \to 0$, as  $\epsilon \to 0$.
   by hypotheses $\bf{H_1} $  and  $\bf{H_2}. $ 
 
 In a similar way, we obtain
 
 \begin{itemize}
  \item $ \displaystyle  \left| \,h_\epsilon^*\nabla_{\Omega_{\epsilon}}h_\epsilon^{*-1} \psi \,(x) 
   \right|_p   
   \leq  B(\epsilon)  |\nabla \psi (x) |_p,    $ 
 \item $ \displaystyle  \left| \,h_\epsilon^*\nabla_{\Omega_{\epsilon}}h_\epsilon^{*-1} \psi \,(x) 
  - \nabla \psi (x) \right|_ p  
    \leq  \eta(\epsilon)  |\nabla \psi (x) |_p,    $ 
   
 \end{itemize}

 It follows that
 
 \begin{align}
     \left| \big\langle \,\big(A_{\epsilon}   -A   \,\big)u\,,\,\psi\,\big\rangle_{ {-1,1}}\right| & 
      \leq 
    \int_{\Omega} \left| \big(\,h_\epsilon^*\nabla_{\Omega_{\epsilon}}h_\epsilon^{*-1}u\,\big)(x) \right|_p \cdot \left| \big(\,h_\epsilon^*\nabla_{\Omega_{\epsilon}}h_\epsilon^{*-1} \psi \,\big)(x) -
     \big( \nabla_{\Omega}  \psi \,\big)(x) \right|_q \, d\,x   \nonumber \\
      &  +    \int_{\Omega}\left| \big(\,h_\epsilon^*\nabla_{\Omega_{\epsilon}}h_\epsilon^{*-1}u -  \nabla_{\Omega} u \,\big)(x) \right|_p \cdot  \left|\big(
    \big( \nabla_{\Omega}  \psi \,\big)(x)  \right|_q \, d\,x   \nonumber \\
  &  +   
 \int_{\Omega}  \left|(h_\epsilon^*\nabla_{\Omega_{\epsilon}}h_\epsilon^{*-1}u)(x) \right|_p \cdot
  \left| \frac{1}{Jh_\epsilon(x) } ( h_\epsilon^* \nabla_{\Omega_{\epsilon}} h_\epsilon^{*-1}  )(x) \, 
 \cdot  \psi (x)  \right|_{q}
  \,dx \nonumber  \\
& \leq  B(\epsilon)  \left[ \int_{\Omega}    |\nabla u (x) |_p^p \, dx \right]^{\frac{1}{p}} \cdot
   \eta(\epsilon)  \left[ \int_{\Omega}    |\nabla \psi (x) |_q^q \, dx \right]^{\frac{1}{q}}  \nonumber  \\
  &  +   \eta(\epsilon)  \left[ \int_{\Omega}    |\nabla u (x) |_p^p \, dx \right]^{\frac{1}{p}} \cdot
    \left[ \int_{\Omega}    |\nabla \psi (x) |_q^q \, dx \right]^{\frac{1}{q}}  \nonumber  \\
    &  +   B(\epsilon) n\cdot \mu(\epsilon)   \left[ \int_{\Omega}    |\nabla u (x) |_p^p \, dx \right]^{\frac{1}{2}} \cdot
    \left[ \int_{\Omega}    | \psi (x) |^q \, dx \right]^{\frac{1}{q}}  \nonumber  \\
    & \leq  \left(  (1+B(\epsilon)) \eta(\epsilon)  + n \beta(\epsilon) ) \cdot \mu(\epsilon)\right) \|u \|_{W{^1,p} (\Omega)}
     \cdot \| \psi \|_{W^{1,q}(\Omega)} \nonumber  \\
& \leq     K(\epsilon) \|u \|_{W{^1,p} (\Omega)}
     \cdot \| \psi \|_{W^{1,q}(\Omega)} \nonumber
      \end{align}

 with  $\displaystyle\lim_{\epsilon \to 0{^+}} K(\epsilon)=0$ (independently of $u$).
 We conclude that 

\begin{align} \label{dif_estimate}
\| \left({A}_{\epsilon} - {A}  \right)u \|_{ {W^{1,q}(\Omega)}'}
& \leq  K(\epsilon)||\,u\,||_{W^{1,p}(\Omega)} \nonumber \\
& \leq \tau(\epsilon)||\,{A} u\,||_{ {W^{q}(\Omega)}'}
\end{align}
 with  $ \displaystyle\lim_{\epsilon \to 0{^+} } \tau(\epsilon)=0$,
 (and $\tau(\epsilon)$ does not depend on $u$).
 \eproof  
 
 \subsection{Existence and continuity of the linear semigroup}
 
 Using  well known facts about the "unperturbed operator"  $A$
  and Theorem 
 \ref{fund_inequality_Lp},
  one can now establish existence and continuity of the linear semigroup, based on the following results:

 \begin{lema} \label{sector}
Suppose  $A $ is a sectorial operator  with  $ \|(\lambda -
A)^{-1}\| \leq  \displaystyle\frac{M}{|\lambda -a |} $ for all $\lambda$ in the
sector $ S_{a,\phi_0} = \{ \lambda  \ | \ \phi_0 \leq
|arg(\lambda-a)| \leq \pi, \lambda \neq a \} $, for some $a\in
\mathbb{R}$ and $0\leq \phi_0<\pi/2$. Suppose also that $B$  is a
linear operator  with  $D(B ) \supset D(A)$ and  $\|Bx - A x\|
\leq \varepsilon \|Ax\| + K \|x\| $, for any  $x \in D(A) $, where
$K$  and  $ \varepsilon $  are positive constants with $
\varepsilon \leq  \displaystyle\frac{1}{4(1+LM)},\, K \leq
\displaystyle\frac{\sqrt{5}}{20M} \frac{\sqrt{2}L - 1}{L^2 - 1} $, for some
$L>1$.

Then  $B$ is also sectorial.  More precisely, if  $b =
\displaystyle\frac{L^2}{L^2 -1} a - \frac{\sqrt{2} L}{L^2 -1} |a| $, $ \phi =
\max \left\{ \phi_0, \displaystyle\frac{\pi}{4}\right\}$ and   $ M' =  2 M \sqrt{5}  $
then
\[
\|(\lambda - B)^{-1}\| \leq  \frac{M'}{|\lambda -b|},
 \]
in the sector  $ S_{b,\phi} = \{ \lambda \ | \ \phi \leq
|arg(\lambda-b)| \leq \pi, \lambda \neq b \} $.
\end{lema}

 {\proof See \cite{PP}, pg 348. \eproof}

\begin{rem} \label{positive} \rm
Observe that $b$ can be made arbitrarily close to $a$ by taking
$L$ sufficiently large. In particular, if $a>0$ then $b>0$.
\end{rem}

\begin{teo} \label{cont_lin_semigroup}
  Suppose that $A$ is as in Lemma \ref{sector}, $\Lambda$ a topological
space and $\{A_\gamma\}_{\gamma \in \Lambda}$
  is a family of operators in $X$ with $A_{\gamma_0} = A$  satisfying the
following conditions:
\begin{enumerate}
  \item  $D(A_{\gamma} ) \supset D(A)$, for all $\gamma\in \Lambda$;
  \item  $\|A_{\gamma}x - A x\|  \leq \epsilon(\gamma) \|Ax\| + K(\gamma)
\|x\| $ for any $x \in D(A)$,
  where $K(\gamma)$ and $\epsilon(\gamma)$ are positive functions with
$\displaystyle\lim_{\gamma \to \gamma_0} \epsilon(\gamma)=0$
  and $\displaystyle\lim_{ \gamma \to \gamma_0} K(\gamma)  = 0 $.
\end{enumerate}
  Then, there exists a neighborhood $V$ of $\gamma_0$ such that
$A_{\gamma}$ is sectorial if $\gamma \in V$ and
 the family of (linear) semigroups $e^{-tA_{\gamma}}$  {satisfies} 
\begin{align} \label{cont_lin_semigroup_eq}
\|  e^{-tA_{\gamma}}-  e^{-tA} \|   \leq        C(\gamma)   e^{-bt}  \nonumber \\
\|  A \left( e^{-tA_{\gamma}}-  e^{-tA} \right) \|   \leq
C(\gamma) \frac{1}{t} e^{-bt} \nonumber \\
\| A^{\alpha} \left( e^{-tA_{\gamma}}-  e^{-tA} \right) \|  \leq
C(\gamma)
\frac{1}{ t^{\alpha}} e^{-bt}, \quad  0< \alpha < 1
\end{align}
for $t >0$, where $b$ is as in Lemma \ref{sector} and
$C(\gamma) \to 0$  as $\gamma \to \gamma_0$.
\end{teo} 

 {\proof See \cite{PP}, pg 349. \eproof}
 
 \begin{teo} \label{family_sec_Lp}
 The  operators  $A_{\epsilon}$  given  by  \eqref{operator_weak_Lp}  in the space $ X= (W^{1,q})'$, with domain 
$W^{1,p}$, $ 1 < p < \infty, \frac{1}{p} +  \frac{1}{p} =1$,    are sectorial operators with sectors and constant in the sectorial inequality independent of $\epsilon$, for  $\epsilon_0$ sufficiently small.   The family  of analytic linear semigroups $e^{- tA_{\epsilon}}$ generated by $A_{\epsilon}$ 
 in the ``base  space''  $X$, satisfies 
 \eqref{cont_lin_semigroup_eq}.  
 \end{teo}
 \proof 
 The first assertion follows from Theorem \ref{sector} and the second from Theorem
 \ref{cont_lin_semigroup},
  \eproof

  \section{The abstract problem in a scale of Banach spaces} 
\label{abstract}
 Our goal in this section is to pose the  problem  (\ref{nonlinBVP})  in a convenient abstract setting.
 
 We proved in Theorem \ref{sectorial_Lp} 
  that, if $\epsilon$ is small,  the operator $A_{\epsilon}$  in $L^p(\Omega)$ defined by
 (\ref{operator_Lp}) with domain given in (\ref{domain_Lp}) is sectorial and, in Theorem \ref{family_sec_Lp} that the same is true for    its extension 
 $\widetilde{A}_{\epsilon}$  to $(W^{1,q})'(\Omega)$.

It is then well-known that the domains $X_{\epsilon}^{\alpha}$ (resp. $\widetilde{X}_\epsilon^{\alpha})$,
 $ \alpha \geq 0$  of the fractional powers of  $A_{\epsilon}$  (resp. $\widetilde{A}_{\epsilon}$) 
are Banach spaces,  $X_{\epsilon}^0 = L^p(\Omega)$, (resp. $\widetilde{X}_\epsilon^0 =  (W^{1,q})'  (\Omega)$),
 $X_{\epsilon}^1 = D(A_{\epsilon}) = W^{2,p} (\Omega) $,  (resp. $\widetilde{X}_\epsilon^1 =  D(\widetilde{A}_{\epsilon}) = 
 W^{1,p} (\Omega) $),
 $X_{\epsilon}^\alpha$,  ( $\widetilde{X}_\epsilon^\alpha$)   is compactly embedded in
  $ X_{\epsilon}^\beta$,   ($  \widetilde{X}_\epsilon^\beta$)  when $ 0 \leq  \alpha < \beta <1 $,
and    $X_{\epsilon}^{\alpha} = W^{p\alpha}$, when $2 \alpha$ is an integer number.


 Since $ X_{\epsilon}^{\frac{1}{2}}= \widetilde{X}_{\epsilon}^1 $, it follows easily that
   $ X_{\epsilon}^{\alpha - \frac{1}{2}  }= \widetilde{X}_{\epsilon}^\alpha  $, for $ \frac{1}{2} 
  \leq \alpha \leq 1 $  and, by an abuse of notation,
 we will still write  $X_{\epsilon}^{\alpha - \frac{1}{2} } $ instead of
  $\widetilde{X}_\epsilon^{\alpha}$,   for  $  0  \leq 
\alpha \leq    \frac{1}{2} $  so we may denote by 
  $ \{ X_{\epsilon}^{\alpha}, \, \,     -\frac{1}{2} \leq  \alpha \leq 1 \} =
 \{ X_{\epsilon}^{\alpha}, \, \,     0 \leq   \alpha \leq 1 \} \cup
 \{ \widetilde{X}_\epsilon^{\alpha}, \, \,   0 \leq   \alpha \leq 1 \},$
 the whole family of fractional power spaces. 
We will denote  simply by  $X^{\alpha}$ the fractional power spaces associated to the unperturbed operator $A$.

 For any  $  \displaystyle -\frac{1}{2} \leq \beta \leq 0$, we may now
define an operator in  these spaces as the restriction of
$\widetilde{A}_{\epsilon}$. We then have the following result
 \begin{teo}\label{sec_scale} For any   $-\frac{1}{2} \leq \beta \leq 0$ and $\epsilon$ sufficiently small, the operator
 $(A_{\epsilon})_{\beta}$ in $ X_{\epsilon}^{\beta}$, obtained by restricting
 $ \widetilde{A}_{\epsilon}$, with domain  $ X_{\epsilon}^{\beta+1}$ is a sectorial operator.
\end{teo}
\proof
 Writing  $\beta =  -\frac{1}{2} + \delta $, for some $ 0 \leq \delta \leq
\frac{1}{2}$, we have $(A_{\epsilon})_{\beta} =
   \widetilde{A}_{\epsilon}^{-\delta} \widetilde{A}_{\epsilon}  \widetilde{A}_{\epsilon}^{\delta}. $ Since $\widetilde{A}_{\epsilon}^{\delta}$ is an isometry from
$X_{\epsilon}^{\beta}$ to $X_{\epsilon}^{ {-\frac{1}{2}}} = (W^{1,q} (\Omega))',$ the result follows easily.  
 {\eproof}

 We
 can now  pose the  problem (\ref{nonlinBVP_fix})  as an abstract problem  in the  scale of Banach spaces $ \{{X}^{\beta}_{\epsilon},
 \, \frac{-1}{2} \leq  \beta  \leq 0  \}  $.

\begin{equation}\label{abstract_scale}
\left\{
\begin{aligned}
u_t  & = -  (A_{\epsilon})_{\beta}u + (H_{\epsilon})_{\beta}u \, , \, t>t_0\,;
\\
u(t_0 & )=u_0 \in X_{\epsilon}^{\eta}\, ,
\end{aligned}
\right.
\end{equation}

where 
\begin{equation} \label{defH}
  (H_{\epsilon})_{\beta} = H(\cdot,\epsilon):=(F_{\epsilon})_{\beta}
  + (G_{\epsilon})_{\beta} :X_{\epsilon} ^{\eta} \to X_{\epsilon} ^{\beta},
  \ \ \epsilon >0 \textrm{ and } 0 \leq \eta \leq \beta +1,
  \end{equation}
\begin{itemize}
\item[(i)] $(F_{\epsilon})_{\beta} = F(\cdot,\epsilon) :X_{\epsilon} ^{\eta} \to X_{\epsilon} ^{\beta}$ is given by
\begin{eqnarray}\label{Fh}
\left\langle F(u,\epsilon)\,,\,\Phi \right\rangle_{\beta \,,\, - \beta} =\displaystyle\int_{\Omega} f(u)\,\Phi\,dx,  \ \ \textrm{ for any } \Phi  \in  (X_{\epsilon}^{  \beta})',
\end{eqnarray}
\item[(ii)] $(G_{\epsilon})_{\beta} = G(\cdot, {\epsilon}) :
 X_{\epsilon} ^{\eta} \to X_{\epsilon} ^{\beta} $ is given by 
\begin{eqnarray}\label{Gh}
\left\langle G(u,\epsilon)\,,\,\Phi \right\rangle_{\beta \,,\,  - \beta } =\displaystyle\int_{\partial \Omega}g(\gamma(u))\,\gamma(\Phi)\left|\frac{J_{\partial \Omega}h_\epsilon}{Jh_\epsilon}\right|\,d\sigma(x)\,, \ \ \textrm{ for any } \Phi \in (X_{\epsilon}^{  \beta})',
\end{eqnarray}
where $\gamma$ is the trace map  {and $J_{\partial\Omega}h_\epsilon$ is the determinant of the Jacobian matrix of the diffeomorphism $h_\epsilon: \partial\Omega \longrightarrow \partial h_\epsilon(\Omega)$}.
\end{itemize}

We will choose $\beta$, small enough in order that $X_{\epsilon}^{\beta+1}$ does not incorporate the boundary conditions, that is, the closure of the subset defined by smooth functions with Neumann boundary condition is the whole space).
It is not difficult to show, integrating by parts, that a regular enough  solution of 
\eqref{abstract_scale}, must satisfy  \eqref{nonlinBVP_fix} 
 (see \cite{COPR} or \cite{OP}).

\section{Local well-posedness}
\label{wellposed}

  In order to prove local  well-posedness for the abstract problem, without assuming growth conditions in 
  the nonlinearities, we want to have two somewhat conflicting requirements for our phase space: we need it to be   continuously embedded in $L^{\infty} $ and we also  do not want it to incorporate the boundary conditions.  To this end, we need to choose $\eta$  and $p$ big enough sot that  the hypotheses of Theorem
  \ref{inclusion} hold   and, on the other hand, we need $\eta$ small enough so that the normal derivative does not have a well defined trace.  To achieve both requirements   we will henceforth 
     assume that 
that 
 \begin{equation} \label{hip_inclusion}
  p  \textrm{ and }  \eta \textrm{ are such the inclusion }  \ \eqref{includ_C}   \
   \textrm{ holds, for some }  \mu \geq 0 \quad \textrm{ and }   \eta < \frac{1}{2}.
 \end{equation} 
  It is easy to check that  \ref{hip_inclusion} holds, for instance, if 
 $p = 2n$, \textrm{ and }
 $    \frac{1}{4} < \eta <  \frac{1}{2}  $. Also,  the last inequality is automatically attended if we choose our 
  base space  $X_{\epsilon}^{\beta} = X_{\epsilon}^{-\frac{1}{2}} = (W^{1,q})' $, where $q $ and $p$ are conjugate exponents, since we must have $\eta- \beta <1 . $

\begin{lema}\label{Flip}
 Suppose that $p $ and $\eta$ are such that  \eqref{hip_inclusion} holds and $f$ is locally Lipschitz.  
Then, the operator 
 $ {(F_{\epsilon})_\eta}   :X_{\epsilon}^{\eta}   {\rightarrow}   X_{\epsilon}^{-\frac{1}{2}}$ given by 
 (\ref{Fh}) is well defined and Lipschitz in bounded sets.
\end{lema}
\proof
 Suppose
 $u \in X_{\epsilon}^{\eta} $.  From \eqref{includ_C} and the hypotheses it follows 
  that $u \in  L^{\infty}(\Omega)  $ and, therefore, if $L_f$  is the Lipschitz constant of $f$ in the interval 
 $ [ -\| u \|_{\infty} ,  \| u \|_{\infty}  ] ]$,
    it follows that
   $ | f(u(x) - f(0) | \leq  L_f | u(x) | $, for  any $x \in \Omega$.
If 
 $\Phi \in (X_{\epsilon}^{-\frac{1}{2}})' =   W^{1,q}$, 

\begin{align*}
\big|\left\langle  {(F_{\epsilon})_\eta}(u)  \,,\,\Phi \right\rangle_{ {\beta, -\beta}}\big| &\leq \displaystyle\int_\Omega |\,f(u)\,|\,|\,\Phi\,|\,dx  \nonumber \\
&\leq L_f \displaystyle\int_\Omega |\,u\,|\,|\,\Phi\,|\,dx  + \displaystyle\int_\Omega |\,f(0)\,|\,|\,\Phi\,|\,dx  \nonumber \\
&\leq L_f   \| u \|_{L^p(\Omega)} \cdot   \| \Phi  \|_{L^q(\Omega)} + \|  f(0)\|_{L^p(\Omega)} \cdot   \| \Phi  \|_{L^q(\Omega)}
\end{align*}

\par  Since  $ W^{1,q} \subset L^q(\Omega)$ and  $X_{\epsilon}^{\eta} \subset L^p(\Omega) $ 
 with stronger norms, we have

\[
\big|\left\langle {(F_{\epsilon})_\eta}(u)\,,\,\Phi \right\rangle_{ { \beta \,,\,-\beta}}\big| 
\leq L_f \|\,u \|_{L^p(\Omega)}||\,\Phi\,||_{W^{1,q}} +  \| f(0)\|_{L^{p} (\Omega)} ||\,\Phi\,||_{W^{1,q}},
\]
so $(F_{\epsilon})_\eta$ is well defined and 

\begin{align}
 \|(F_{\epsilon})_{\eta}(u)\|_{(W^{1,q}   )' }  & \leq  L_f  \|u\ \|_{L^p(\Omega)}  + \| f(0)\|_{L^{p} (\Omega)}  \label{normF_Lp} \\
  & \leq  L_f  \|u\ \|_{X_{\epsilon}^\eta}  + \| f(0)\|_{L^{p} (\Omega)}  \label{normF}
\end{align}
 where $L_f$  is the Lipschitz constant of $f$ in the interval 
 $ [ -\| u \|_{\infty} ,  \| u \|_{\infty}  ] ]$

 Alternatively, if $M_f = M_f(u): = \sup \{|f(x)| \ x \in
   [ -\| u \|_{\infty} ,  \| u \|_{\infty}  ] ] \}$, 
    it follows that
 
\begin{align*}
  \big|\left\langle  {(F_{\epsilon})_\eta}(u)  \,,\,\Phi \right\rangle_{ {\beta, -\beta}}\big| &\leq \displaystyle\int_\Omega |\,f(u)\,|\,|\,\Phi\,|\,dx
                                                                                                 \nonumber \\
                                                                                               &\leq  M_f  \displaystyle\int_\Omega  \,|\,\Phi\,|\,dx
                                                                                                 \nonumber \\
  & \leq  M_f  |\Omega|^{\frac{1}{p} }  \,|\,\Phi\,|_{L^q(\Omega)}   \\
   & \leq  M_f  |\Omega|^{\frac{1}{p} }  \,|\,\Phi\,|_{W^{1,q}(\Omega)}   \\
\end{align*}

Thus 
\begin{align}
  \|(F_{\epsilon})_{\eta}(u)\|_{(W^{1,q}   )' }   & \leq  M_f   |\Omega|^{\frac{1}{p} } 
    \label{normF_infty}
\end{align}

  Suppose now that If  $u_1, u_2$ belong to a bounded set $B \in X_{\epsilon}^{ {\eta}}$. 
   From \eqref{includ_C} and the hypotheses it follows now that $u_1, u_2$ belong to a ball of radius
    $R = \sup_{u\in B} \|u\|_{\infty}$
   in $L^{\infty}(\Omega)  $ and, therefore, if L is the Lipschitz constant of $f$ in the interval $[-R,R]$,
   we have 
   $ | f(u_1(x)) - f(u_2(x))) | \leq L | u_1(x) - u_2(x) | $, for any $x \in \Omega$.
If 
 $\Phi \in (X_{\epsilon}^{-\frac{1}{2}})' =   W^{1,q}$, we obtain

\begin{align*} 
\left|\left\langle{(F_{\epsilon})_\eta}(u_1) - {(F_{\epsilon})_\eta}(u_2) \,,\,\Phi \right\rangle_{\beta, -\beta}\right| & = 
\left|\displaystyle\int_{\Omega} [f(u_1)-f(u_2)]\,\Phi\,dx\,\right| \\
&\leq  \displaystyle\int_{\Omega} L  \,|\,u_1-u_2\,|\,|\,\Phi\,|\,dx  \\
&\leq   L_f \| u_1 - u_2 \|_{L^p(\Omega)} \cdot   \| \Phi  \|_{L^q(\Omega)} \\
&\leq   L_f \| u_1 - u_2 \|_{X_{\epsilon}^\eta}  \cdot   \| \Phi  \|_{W^{1,q}(\Omega)} 
\end{align*} 

Thus

\begin{align} 
\|(F_{\epsilon})_{\eta}(u_1) -  (F_{\epsilon})_{\eta}(u_2)\|_{(W^{1,q})' } & \leq 
 L_f \| u_1 - u_2 \|_{L^p(\Omega)}  \label{LipF_Lp}\\ 
 & \leq  L_f \| u_1 - u_2 \|_{X_{\epsilon}^\eta}. \label{LipF}
\end{align} 

 This concludes the proof.

\eproof

\begin{lema}\label{Gbem}
 Suppose that $p $ and $\eta$ are such that  \eqref{hip_inclusion} holds and $g$ is locally Lipschitz.  
 Then, if $\epsilon_0$ is sufficiently small,  the operator $(G_{\epsilon})_{\eta}  {=G} :X_{\epsilon}^{\eta}   {\rightarrow} (W^{1,q})' $
given by (\ref{Gh})  is well defined, for $0 \leq \epsilon <\epsilon_0$
and  bounded in bounded sets.
\end{lema} 
\proof 
 Suppose
 $u \in X_{\epsilon}^{\eta} $.  From \eqref{includ_C} and the hypotheses it follows  that $u \in  L^{\infty}(\Omega)  $ and, therefore, if $L_g$  is the Lipschitz constant of $g$ in the interval  
 $ [ -\| u \|_{\infty}, \| u \|_{\infty}    ]$,
    it follows that  
   $ | g(\gamma(u)(x) - g(0) | \leq L_g | \gamma(u)(x) | $, for  any $x \in \partial \Omega$.

 If $u \in X_{\epsilon}^{\eta}$ and   $\Phi \in (X_{\epsilon}^{-\frac{1}{2}})' =   W^{1,q}$, we have

\begin{align*}
  \big|\left\langle G(u,\epsilon)\,,\,\Phi \right\rangle_{\beta \,,\, -\beta }\big| &\leq\displaystyle\int_{\partial \Omega} |\,g(\gamma(u))\,|\,|\,\gamma(\Phi)\,
  |\left|\frac{J_{\partial \Omega}h_\epsilon}{Jh_\epsilon}\right|\,d\sigma(x) \nonumber  \\
  &\leq ||\mu||_\infty \displaystyle\int_{\partial \Omega}  L_g |\gamma(u)| |\gamma(\Phi)| +
    |g(0)|\,|\gamma(\Phi)| \, d\sigma(x) \nonumber  \\
    &\leq ||\mu||_\infty  \left( L_g  \|\gamma(u) \|_{L^p(\partial \Omega)} \cdot 
    \|\gamma(\Phi)\|_{L^q(\partial \Omega)}
     +  \| g(0) \|_{L^p(\partial \Omega)} \cdot 
    \|\gamma(\Phi)\|_{L^q(\partial \Omega)} \right) 
\end{align*}

where  $\mu(x, \epsilon) = \left|\frac{J_{\partial \Omega}h_\epsilon}{Jh_\epsilon}\right|$, and  $\|\mu\|_\infty = \sup \left\{ |\mu(x, \epsilon)| \, | \, 
 x\in \partial \Omega, \,  0 \leq \epsilon
 \leq  \epsilon_0 \right\}$                    is finite by
  hypothesis $ ({\bf H_1})$.
  
By  Theorem \ref{trace}, there

$||\gamma(\Phi)||_{L^q(\partial \Omega)} \leq  {K}_1\,||\Phi||_{W^{1,q}(\Omega)} \,,\
||\gamma(u)||_{L^p(\partial \Omega)} \leq {K}_2 \,||u||_{X_{\epsilon}^{\eta}}\, ,
$
for some constants  ${K}_1$,   ${K}_2$.

   Thus

 \begin{align*}
  \big|\left\langle G(u,\epsilon)\,,\,\Phi \right\rangle_{\beta \,,\, -\beta }\big| 
  &\leq ||\mu||_\infty \left(  L_g K_1   \,||\gamma(u)||_{L^p(\partial \Omega) }
    \|\Phi |_{W^{1,q} (\Omega)} 
     + K_1 \| g(0) \|_{L^p(\partial \Omega)} \cdot 
    \|\Phi \|_{W^{1,q}(\Omega) } \right)  
\end{align*}

proving that  $(G_{\epsilon})_{\beta} $ is well defined and 
 \begin{align}  
 \|G(u,\epsilon) \|_{(W^{1,q}(\Omega))'   } &  \leq 
||\mu||_\infty  \left( L_g  K_1   \,||\gamma(u)||_{ L^p(\partial \Omega) } 
     + K_1 \| g(0) \|_{L^p(\partial \Omega)} \right)   \label{normG_Lp}  \\
      &  \leq 
||\mu||_\infty  \left( L_g  K_1 K_2  \,||u||_{X_{\epsilon}^{\eta}} 
     + K_1 \| g(0) \|_{L^p(\partial \Omega)} \right)  \label{normG} 
  \end{align}

 Alternatively, if $M_g= M_g (u): = \sup \{|g(x)| \ x \in
   [ -\| u \|_{\infty} ,  \| u \|_{\infty}  ] ] \}$, 
    it follows that

    \begin{align*}
  \big|\left\langle G(u,\epsilon)\,,\,\Phi \right\rangle_{\beta \,,\, -\beta }\big| &\leq\displaystyle\int_{\partial \Omega} |\,g(\gamma(u))\,|\,|\,\gamma(\Phi)\,
  |\left|\frac{J_{\partial \Omega}h_\epsilon}{Jh_\epsilon}\right|\,d\sigma(x) \nonumber  \\
                                                                                    &\leq ||\mu||_\infty M_g \displaystyle\int_{\partial \Omega}  \,|\gamma(\Phi)| \, d\sigma(x) \nonumber  \\
                                                                                    &\leq ||\mu||_\infty M_g |{\partial \Omega}|^{\frac{1}{p} }  \, \|\gamma(\Phi)\|_{L^q( \partial \Omega)} \nonumber  \\
       &\leq ||\mu||_\infty M_g |{\partial \Omega}|^{\frac{1}{p} }  \,K_1 \|\Phi) \|_{W^{1,q}( \Omega)} \nonumber  \\
       \end{align*}

       Thus

        \begin{align}
 \|G(u,\epsilon) \|_{(W^{1,q}(\Omega))'   }
       &\leq ||\mu||_\infty M_g |{\partial \Omega}|^{\frac{1}{p} }  \,K_1  \label{normG_infty}
       \end{align}

\eproof

\begin{lema}\label{Glip}
Suppose the same hypotheses of Lemma  \ref{Gbem} hold. 
Then
 the operator
 $G(u,\epsilon)   {=G(u)} :X_{\epsilon}^{\eta} \times [0, \epsilon_0] \to (W^{1,q})'$ given
 by  (\ref{Gh})
 is uniformly continuous in $\epsilon$, for $u$ in bounded sets of $X_{\epsilon}^{\eta}$
   and
  locally Lipschitz continuous in $u$, uniformly in  
  $\epsilon$. 
\end{lema}
\proof
 We first show that $(G_{\epsilon})_{\beta} $  is locally Lipschitz continuous in  $u \in X_{\epsilon}^{\eta}$.  
 
Suppose  that   $u_1, u_2$ belong to a bounded set $B \in X_{\epsilon}^{ {\eta}}$. 
   From \eqref{includ_C}, the Trace Theorem and the hypotheses, it follows now that $\gamma(u_1), \gamma(u_2)$ belong to a ball of some radius $R$
   in $L^{\infty}(\partial \Omega)  $ and, therefore, if $L_g$  is the Lipschitz constant of $g$ in the interval $[-R,R]$,
   we have 
   $ | g(\gamma(u_1)(x)) - g(\gamma(u_2)(x))) | \leq L_g | \gamma(u_1)(x) - \gamma(u_2)(x) | $, for any $x \in \partial \Omega$.
If 
 $\Phi \in     (W^{1,q})'$ and   $ \epsilon \in [0, \epsilon_0]$, we obtain

 Then

\begin{align*}
\left|\left\langle G(u_1,\epsilon) - G(u_2,\epsilon),\Phi \right\rangle_{\beta,-\beta}\right |
 &\leq  \displaystyle\int_{\partial\Omega} |\,g(\gamma(u_1))-g(\gamma(u_2))\,|\,\big|\,\gamma(\Phi)\,\big|\,\left|\frac{J_{\partial\Omega} h_\epsilon}{Jh_\epsilon}\right|\,d \sigma (x) \\
&\leq \displaystyle\int_{\partial\Omega} L_g  |\gamma(u_1)-\gamma(u_2)|\left|\gamma(\Phi)\right|\left|\frac{J_{\partial\Omega} h_\epsilon}{Jh_\epsilon}\right|d\sigma (x) \\
&\leq  
 L_g ||\mu||_\infty  \displaystyle\int_{\partial\Omega}   |\gamma(u_1)-\gamma(u_2)|\left|\gamma(\Phi)\ \right|d\sigma (x) \\
 &\leq   L_g ||\mu||_\infty   \| \gamma(u_1) - \gamma( u_2 ) \|_{L^p(\partial \Omega)} \cdot   \| \gamma(\Phi )  \|_{L^q(\partial \Omega)} \\
 &\leq   L_g   ||\mu||_\infty  K_1 K_2 \| u_1  - u_2 \|_{X_{\epsilon}^{\eta}} \cdot   \| \Phi  \|_{W^{1,q}(\Omega)}, 
\end{align*}
where $K_1, K_2$ are the norms of the trace mappings, given by Theorem \ref{trace}. 
  Therefore,
  \begin{align} 
\| G(u_1,\epsilon) - G(u_2,\epsilon) \|_{(W^{1,q})' }
 &\leq   L_g  ||\mu||_\infty  K_1    \| \gamma(u_1) - \gamma( u_2 ) \|_{L^p(\partial \Omega)}  
 \label{LipG_Lp} \\
 &\leq   L_g  ||\mu||_\infty  K_1 K_2 \| u_1  - u_2 \|_{X_{\epsilon}^{\eta}} \label{LipG}
  \end{align}
  so  $(G_{\epsilon})_{\beta} $ is locally Lipschitz in $u$.

\par Now, if  $u \in X_{\epsilon}^{\eta}$,  $\Phi \in     (W^{1,q})'$   and
 $ \epsilon_1, \epsilon_2 \in [0, \epsilon_0]$, we have
 \begin{align*}
\big|\langle G(u,{\epsilon_1})- & G(u,{\epsilon_2}),\Phi\rangle_{\beta,-\beta}\big| \leq
\displaystyle\int_{\partial\Omega}|\,\gamma(g(u))\,|\,|\,\gamma(\Phi)\,|\left|\,\left(\,\left|\frac{J_{\partial \Omega}h_{\epsilon_1}}{Jh_{\epsilon_1}}\right|-\left|\frac{J_{\partial \Omega}h_{\epsilon_2}}{Jh_{\epsilon_ 2}}\right|\,\right)\,\right|\,d\sigma(x)  \\
&\leq
\|\mu_{\epsilon_1} - \mu_{\epsilon_2} \|_{\infty}
\displaystyle\int_{\partial\Omega} |\,g(\gamma(u))\,|\,|\,\gamma(\Phi)\,|\, d\sigma(x) \\
&\leq
\|\mu_{\epsilon_1} - \mu_{\epsilon_2} \|_{\infty}
\displaystyle\int_{\partial\Omega} \left( L_{g} |\,\gamma(u)\, | +  |\, g(0)\, | \right)  \,|\,\gamma(\Phi)\,|\, d\sigma(x) \\
 &\leq  \|\mu_{\epsilon_1} - \mu_{\epsilon_2} \|_{\infty} \left(  L_g  \|\gamma(u) \|_{L^p(\partial \Omega)} \cdot 
    \|\gamma(\Phi)\|_{L^q(\partial \Omega)}
     +  \| g(0) \|_{L^p(\partial \Omega)} \cdot 
    \|\gamma(\Phi)\|_{L^q(\partial \Omega)} \right) \\
    & \leq 
      \|\mu_{\epsilon_1} - \mu_{\epsilon_2} \|_{\infty} \left(  L_g K_1 K_2  \,||u||_{X^{\eta}}
    \|\Phi |_{W^{1,q} (\Omega)} 
     + K_1 \| g(0) \|_{L^p(\partial \Omega)} \cdot 
    \|\Phi \|_{W^{1,q}(\Omega) } \right)  
\end{align*}

\noindent  where $\|\mu_{\epsilon_1}- \mu_{\epsilon_2} \|_\infty =
 \sup \left\{ \left|\frac{J_{\partial \Omega}h_{\epsilon_1}}{Jh_{\epsilon_1}}\right|-\left|\frac{J_{\partial \Omega}h_{\epsilon_2}}{Jh_{\epsilon_ 2}}\right|  \, | \,
 x\in \partial \Omega, \, 
   \right\} \to 0 $ as $ |\epsilon_1 - \epsilon_2| \to 0 $,
 by hypothesis $ ({\bf H_1})$ and $K_1, K_2$ are trace constants given by  Theorem \ref{trace}.
 
 It follows that  
 \begin{equation}  \label{unif_contG} 
 \|G(u,\epsilon_1) -  G(u,\epsilon_2)\|_{(W^{1,q}(\Omega))'   } \leq 
  \|\mu_{\epsilon_1} - \mu_{\epsilon_2} \|_{\infty} \left(  L_g K_1 K_2  \,||u||_{X^{\eta}} 
     + K_1 \| g(0) \|_{L^p(\partial \Omega)} \right)
 \end{equation}

  Alternatively, if $M_g= M_g (u): = \sup \{|g(x)| \ x \in
   [ -\| u \|_{\infty} ,  \| u \|_{\infty}  ] ] \}$,
   
 \begin{align*}
\big|\langle G(u,{\epsilon_1})- & G(u,{\epsilon_2}),\Phi\rangle_{\beta,-\beta}\big| \leq
\displaystyle\int_{\partial\Omega}|\,\gamma(g(u))\,|\,|\,\gamma(\Phi)\,|\left|\,\left(\,\left|\frac{J_{\partial \Omega}h_{\epsilon_1}}{Jh_{\epsilon_1}}\right|-\left|\frac{J_{\partial \Omega}h_{\epsilon_2}}{Jh_{\epsilon_ 2}}\right|\,\right)\,\right|\,d\sigma(x)  \\
&\leq
\|\mu_{\epsilon_1} - \mu_{\epsilon_2} \|_{\infty}
\displaystyle\int_{\partial\Omega} |\,g(\gamma(u))\,|\,|\,\gamma(\Phi)\,|\, d\sigma(x) \\
&\leq
\|\mu_{\epsilon_1} - \mu_{\epsilon_2} \|_{\infty}  M_g
\displaystyle\int_{\partial\Omega}   \,\gamma(\Phi)\,|\, d\sigma(x) \\
&\leq
\|\mu_{\epsilon_1} - \mu_{\epsilon_2} \|_{\infty}  M_g |\partial \Omega|^{\frac{1}{p}  }
 \|  \,\gamma(\Phi) \|_{L^p(\partial \Omega)}  \\
&\leq
\|\mu_{\epsilon_1} - \mu_{\epsilon_2} \|_{\infty}  M_g |\partial \Omega|^{\frac{1}{p}  }
 K_1 \|  \,\Phi \|_{W^{1,q}}(\Omega) 
\end{align*}

 It follows that  
 \begin{equation}  \label{unif_contG_infty} 
 \|G(u,\epsilon_1) -  G(u,\epsilon_2)\|_{(W^{1,q}(\Omega))'   } \leq 
 \|\mu_{\epsilon_1} - \mu_{\epsilon_2} \|_{\infty}  M_g |\partial \Omega|^{\frac{1}{p}  } K_1
 \end{equation}

 \eproof
 
\begin{cor} \label{propH} Suppose the hypotheses of Lemmas \ref{Flip} and \ref{Gbem} hold.  Then the map
$(H (u,\epsilon))_{\eta} := ( F_{\epsilon} (u))_{\eta} + (G (u,\epsilon))_{\eta}    :X_{\epsilon}^{\eta} \times [0, \epsilon_0] \to (W^{1,q})'$ is well defined, bounded in bounded sets uniformly in $\epsilon$,  uniformly continuous in $\epsilon$ for $u$ in bounded sets of $X_{\epsilon}^{\eta}$
   and
  locally Lipschitz continuous in $u$  uniformly in  
  $\epsilon$. 
\end{cor} 
\proof
 
 From \ref{normF_Lp}, \ref{normF}, \ref{normG_Lp} and  \ref{normG}, we obtain  

 \begin{align}
\|(H_{\epsilon})_{\eta}(u)\|_{(W^{1,q}   )' }  & \leq   L_f  \|u  \|_{L^p(\Omega)} +
 L_g K_1  \|\gamma(u)  \|_{L^p(\partial \Omega)}
   + \| f(0)\|_{L^{p} (\Omega)} 
  +  K_1 \| g(0) \|_{L^p(\partial \Omega)}   \label{normH_Lp} \\
  & \leq  \left(   L_f   +
 L_g K_1  K_2  \right)  ||u||_{X_{\epsilon}^{\eta}}
   + \| f(0)\|_{L^{p} (\Omega)} 
  +  K_1 \| g(0) \|_{L^p(\partial \Omega)},   \label{normH} 
  \end{align} 
 where $L_f$ and $L_g$  are Lipschitz constants of $f$ and $g$  in the interval 
   $ [ -\| u \|_{\infty} ,  \| u \|_{\infty}  ] ]$, 
 
 Alternatively, if $M_f = M_f(u): = \sup \{|f(x)| \ x \in
   [ -\| u \|_{\infty} ,  \| u \|_{\infty}  ] ] \}$, 
   $M_g= M_g (u): = \sup \{|g(x)| \ x \in
   [ -\| u \|_{\infty} ,  \| u \|_{\infty}  ] ] \}$, we obtain from \eqref{normF_infty} and 
   \eqref{unif_contG_infty}

  \begin{align}
\|(H_{\epsilon})_{\eta}(u)\|_{(W^{1,q}   )' }  & \leq  M_f   |\Omega|^{\frac{1}{p} } +
 ||\mu||_\infty M_g |{\partial \Omega}|^{\frac{1}{p} }  \,K_1  \label{normH_infty}
  \end{align}

\vspace{3mm}

  From \eqref{LipF_Lp}, \eqref{LipF}, \eqref{LipG_Lp} and  \eqref{LipG},

\begin{align} 
\| H(u_1,\epsilon) - H (u_2,\epsilon) \|_{(W^{1,q})' }
  &\leq   L_g  ||\mu||_\infty  K_1    \| \gamma(u_1) -
    \gamma( u_2 ) \|_{L^p(\partial \Omega)}  +   L_f \| u_1 - u_2 \|_{L^p(\Omega)}  
 \label{LipH_Lp} \\
  &\leq  \left( L_g  ||\mu||_\infty  K_1 K_2 +  L_f \right)
    \| u_1  - u_2 \|_{X_{\epsilon}^{\eta}} \label{LipH}
  \end{align}

\vspace{3mm}

  From \eqref{unif_contG}
  \begin{equation}  \label{unif_contH} 
 \|H(u,\epsilon_1) -  H(u,\epsilon_2)\|_{(W^{1,q}(\Omega))'   } \leq 
  \|\mu_{\epsilon_1} - \mu_{\epsilon_2} \|_{\infty} \left(  L_g K_1 K_2  \,||u||_{X^{\eta}} 
     + K_1 \| g(0) \|_{L^p(\partial \Omega)} \right)
 \end{equation} 
 
 Alternatively, from \eqref{unif_contG_infty}
\begin{equation}  \label{unif_contH_infty} 
 \|H(u,\epsilon_1) -  H(u,\epsilon_2)\|_{(W^{1,q}(\Omega))'   } \leq 
 \|\mu_{\epsilon_1} - \mu_{\epsilon_2} \|_{\infty}  M_g |\partial \Omega|^{\frac{1}{p}  } K_1
 \end{equation}     
 
  In the estimates above  ${K}_1$ and  ${K}_2$ are the norms of  the trace  mappings (see Theorem 
  \ref{trace}).
 \eproof

\begin{teo}\label{loc_exist}
  Suppose  the hypotheses of Corollary \ref{propH}  hold.  
 Then, for any $(t_0, u_0) \in \R \times X_{\epsilon}^\eta$, the  problem (\ref{abstract_scale}) 
 has a unique solution $u(t,t_0, u_0,\epsilon)$  with initial value $u(t_0) = u_0$.
\end{teo}
\proof
 From Theorem \ref{sec_scale}  it follows that $ {(}A_{ {\epsilon}} {)_\beta}$ is a sectorial operator  in
 $(W^{1,q})'$, with domain  $X_{\epsilon}^{\frac{1}{2}} =  W^{1,p} $, if $\epsilon $ is small enough. 
 The result follows then from Corollary \ref{propH}  and   results in \cite{He2}
 and \cite{PP}.

 \eproof

\section{Global existence and boundedness of the semigroup} 
\label{global_exist}

  We will use the notation $T_{\epsilon}(t) u_0$ for the  (local) solution of the problem \eqref{abstract_scale} given by
   Theorem \ref{loc_exist}, with initial condition $u_0$ in some fractional power space of $A_{\epsilon}$. 
 We now want to show that these solutions are globally defined if an additional (dissipative)  hypotheses on $f$ and $g$ is assumed.
 Here are these hypotheses:
 
%
 
There exist constants  $c_0$ and  $d_0$ such that
\begin{equation} \label{dissipative}
\displaystyle\limsup_{|\,u\,| \to \infty }\frac{f(u)}{u} \leq c_0\,, \quad
\displaystyle\limsup_{|\,u\,| \to \infty }\frac{g(u)}{u} \leq d_0 
\end{equation}
and  the first eigenvalue $\mu_1(\epsilon)$ of the problem
\begin{equation} \label{compet}
\left\{
\begin{array}{lll}
- h_{\epsilon}^{*}  \Delta_{\Omega_\epsilon} {h_{\epsilon}^{*}}^{-1} \Delta u + (a - c_0)u = \mu u \,\,\mbox{em} \,\,\Omega   \\
\displaystyle h_{\epsilon}^{*} \frac{\partial u}{\partial N_\Omega}{h_{\epsilon}^{*}}^{-1} =d_0\,u \,\, \mbox{em} \,\, \partial\Omega
\end{array}
\right.
\end{equation}
is positive for $\epsilon$ sufficiently small.

\begin{rem} \label{dissip_perturbed}
 Observe that if   the hypothesis  (\ref{compet})  hold
 for $\epsilon = 0$, then this also true for $\epsilon$ small   since  the eigenvalues change continuously with $\epsilon$ by (\ref{dif_estimate}).
\end{rem}
  
 \begin{rem}
The arguments bellow are a slight modification of the ones in 
{\cite{OP}}, but we include  them  here for the sake of completeness. Similar arguments were used in \cite{ACB} in a somewhat different setting.
\end{rem}

 In order to   use comparison results, we start 
by defining the concepts of sub- and super-solutions.

\begin{defn}
Suppose $\Omega$  is a $\mathcal{C}^{1,\alpha}$, domain for some $\alpha \in (0,1)$, $L$ is a uniformly elliptic  second order differential operator in $\overline{\Omega}$,
  $u_0\in \mathcal{C }^{\alpha}(\Omega)$, $T>0$ and
$\bar{u}:\Omega\subset\R^\to\R^n$ (${\underline u}$
respectively)  a function which is continuous in $[0,T]\times\bar\Omega$, 
continuously differentiable in $t$ and twice continuously differentiable in
$x$
for 
$(t,x)\in (0,T]\times\Omega$. Then $\overline{u}$ (${ respectively, \underline u}$)
is a super-solution (sub-solution) of the problem

\begin{equation}
 \left\{
\begin{aligned}
\displaystyle
  u_t & = Lu +f(u), \qquad \hbox{in}  \quad (0,T]\times \Omega, \\  
 \frac{\partial u}{\partial N} & = g(u), \qquad \hbox{on}\quad \partial \Omega 
\\
  u(0) & =u_0.
\label{subsup}
\end{aligned}
 \right.
\end{equation}
if it satisfies

\begin{equation}
 \left\{
\begin{aligned}
\displaystyle
  u_t & \geq  Lu +f(u) , \qquad \hbox{in}  \quad (0,T]\times \Omega, \\  
 \frac{\partial u}{\partial N} & \geq  g(u), \qquad \hbox{on}\quad \partial \Omega
\\
  u(0) &  \geq u_0.
\label{super}
\end{aligned}
 \right.
\end{equation}

( and respectively with the $\geq$ sign replaced by the $\leq$ sign).
\end{defn}

 A basic result for our arguments is the following
\begin{teo}(Pao \cite{Pao})
\label{pao}

If $f$ is locally Lipschitz and ${\bar u}$ and ${\underline u}$ are
respectively 
a super  and sub-solution of the problem (\ref{subsup}), satisfying
\[
{\underline u}\leq {\bar u},{\hbox{ in }}\Omega\times(0,T),
\]
then there exists a solution $u$ of (\ref{subsup}) such that
\[
{\underline u}\leq u\leq {\bar u},{\hbox{ in }}\Omega\times(0,T).
\]
\end{teo}

\vspace{3mm}

Let now  $\varphi_{\epsilon}$ be the first positive normalized 
eigenfunction of  \eqref{compet} and
$\displaystyle m_{\epsilon} =\min_{x\in{\bar\Omega}}\varphi_{\epsilon} (x)$. We know that $m_{\epsilon}>0$. 
For each $\theta > 0 \in\R $, define
\[
\Sigma_\theta^{\epsilon}=\displaystyle\left\{u\in X^{\eta}_{\epsilon} : |u(x)|\leq \theta\varphi_{\epsilon} (x),
{\hbox{ for all }}x\in{\bar\Omega}\right\}.
\]

 From the dissipative hypothesis \eqref{dissipative} on $f$ and $g$, we know that there
exists
$\xi\in\R$, such that
\[
\displaystyle
 \frac{f(s)}{s}\leq c_0{\hbox{ and }} \frac{g(s)}{s}\leq d_0,
\]
for all $s$ with $|s|\geq\xi$.
 To simplify the notation, we take the $\epsilon = 0$,  in the proofs below, since the argument is the same for any $\epsilon$ such that  \eqref{compet} is true (see Remark \ref{dissip_perturbed}).
\begin{lema} 
  \label{posinv}  Suppose, in addition to
  the hypotheses of Theorem \ref{loc_exist}, that \eqref{dissipative}  and
   \eqref{compet} hold. 
   Then, if $\theta m_{\epsilon} \geq\xi$  and $\epsilon$ is small enough, the set
   $\Sigma_\theta^{\epsilon} $ is a positively invariant set for $T(t)$.
\end{lema}

\medskip

\proof

Let 
\[
\begin{array}{l}
\displaystyle \Sigma^1_\theta =\{u\in X^{\eta} : u(x)\leq \theta\varphi(x),{\hbox{
for all }}
x\in{\bar\Omega}\}\\
\\
\displaystyle \Sigma^2_\theta =\{u\in X^{\eta} : u(x)\geq -\theta\varphi(x),{\hbox{
for all }}
x\in{\bar\Omega}\}
\end{array}
\]

Since $\Sigma_\theta=\Sigma^1_\theta\cap\Sigma^2_\theta$ it is enough to show 
that $\Sigma^1_\theta$ and $\Sigma^2_\theta$ are positively invariant.

Let $u_0\in\Sigma^1_\theta$, and suppose, for contradiction, that there 
exists $t_0\in[0,t_{\max}[$ and $x_0\in{\bar\Omega}$ such that
\[
T(t_0)u_0(x_0) > \theta\varphi(x_0).
\]
Consider ${\bar v}(t)=e^{-\mu(t-t_0)}\theta\varphi$, where $\mu$ is the 
eigenvalue associated with $\varphi$. We have that

\[ \left\{
\begin{aligned}
  \frac{\partial{\bar v}}{\partial t} &=   \left(  \Delta  {\bar v}- 
   a{\bar v}+c_0{\bar v}   \right)
  \geq 
  \Delta  {\bar v}
  - a{\bar v}+ f({\bar v}) 
\\
\frac{\partial {\bar v}}{ \partial N} & =  d_0{\bar v} \geq g({\bar v}),
\end{aligned}
\right.
\]
for all $t\in]0,t_0]$.

Thus ${\bar v}$ is a super-solution for the problem (\ref{nonlinBVP_fix}). It follows 
from Theorem~{\ref{pao}} that 
\[
T(t)u_0\leq {\bar v}(t),{\hbox{ in }}{\bar\Omega}{\hbox{ for all
}}t\in[0,t_0[.
\]
In particular, $T(t_0)u_0(x_0)\leq\theta\varphi(x_0)$ and we reach a
contradiction.

To prove that $\Sigma^2_\theta$ is positively invariant we proceed in a
similar
way, 
using now that ${\underline v}=-{\bar v}$ is a sub-solution for the problem 
(\ref{nonlinBVP_fix}).

\eproof

\medskip

\begin{lema}
\label{bound_theta}
Suppose the hypotheses of Lemma \ref{posinv} hold. 
If $\theta m_{\epsilon} \geq \xi$,  and $ \eta \leq \alpha <  \frac{1}{2} $, there exists a constant $R= R(\theta, \eta)$, and $T>0$  independent of $\epsilon$,  such that the orbit of any
  bounded subset V  of $X^{\eta}_{\epsilon} \cap \Sigma_\theta^{\epsilon}$  under $T_{\epsilon}(t)$ is  in the ball of radius 
   $R$ of   $X^{\alpha}_{\epsilon}$, for $ t > T $. In particular, the solutions with initial condition in $X^{\eta}_{\epsilon} \cap \Sigma_\theta$  are globally defined. 
\end{lema}
\medskip

\noindent{\bf Proof:}

\medskip

 Lemma~\ref{posinv} implies that 
$T_{\epsilon}(t)u_0\in\Sigma_\theta^{\epsilon}$, for all $t\in[0,t_{\max}[$ so 
\[
\|T_{\epsilon}(t)u_0\|_\infty\leq\theta\|\varphi\|_\infty.
\]

Applying the variation of constants formula, we obtain (see \cite{He2})
\[
\displaystyle
\|T(t)u_0\|_{X^{\alpha}}\leq M  t^{-(\alpha - \eta  ) } e^{-\delta
t}\|u_0\|_{X^{\eta}}+M\int_0^t(t-s)^{-(\alpha +\frac{1}{2})}
 e^{-\delta(t-s)}\|(H_{\epsilon})_{\eta}(T(s)u_0)\|_{X^{-\frac{1}{2}    }    }     \,ds,
\]
where the   $M,\delta >0$ are constants depending only on the decay of the linear semigroup 
$e^{A_{\epsilon}t}$, and can be chosen independently of $\epsilon$.
By \eqref{normH_infty}

 \begin{align*}
\|(H_{\epsilon})_{\eta}(T(s)u_0)\|_{X^{-\frac{1}{2}  }}  &  \leq  M_f   |\Omega|^{\frac{1}{p} } +
 ||\mu||_\infty M_g |{\partial \Omega}|^{\frac{1}{p} }  \,K_1  \\
 \end{align*}
where now  
$M_f = M_f(u): = \sup \{|f(x)| \ x \in I \}$, 
   $M_g= M_g (u): = \sup \{|g(x)| \ x \in I \}$
with  $ [ - \theta \| \phi_{\epsilon}\|_{\infty} , \phi_{\epsilon}\|_{\infty}   ] \subset I $, for all
 $\epsilon$ sufficiently small. 
Thus, writing $K=  M_f   |\Omega|^{\frac{1}{p} } +
 ||\mu||_\infty M_g |{\partial \Omega}|^{\frac{1}{p} }  \,K_1$, we obtain

\begin{align*}
\displaystyle
\|T_{\epsilon}(t)u_0\|_{X^{\alpha}} & \leq  M    t^{-(\alpha - \eta  ) } e^{-\delta
t}\|u_0\|_{X^{\eta}}+ K M\int_0^t(t-s)^{-(\alpha +\frac{1}{2})}
 e^{-\delta(t-s)}        \,ds,  \ \\
  & \leq 
 M    t^{-(\alpha - \eta  ) } e^{-\delta
t}\|u_0\|_{X^{\eta}}+ 
KM \frac{ \Gamma( \frac{1}{2} -\alpha   ) }{\delta^{\frac{1}{2}- \alpha} }.    
  \end{align*}

for all $t\in[0,t_{max}[$.

Therefore $\|T_{\epsilon}(t)u_0\|_{X^{\alpha}}$ is bounded by a constant for any $t > 0$.   Since
$ X^{\alpha}$ is compactly embedded in $ X^{\eta }$, if $\alpha > \eta$, it follows that the solution is globally defined.  Also, if $T$ is such that 
 $   t^{-(\alpha - \eta  ) } e^{-\delta
t}\|u_0\|_{X^{\eta}} \leq  
K  \frac{ \Gamma( \frac{1}{2} -\alpha   ) }{\delta^{\frac{1}{2}- \alpha} } $, then 
$\|T_{\epsilon}(t)u_0\|_X^{\alpha}$ belongs  to the ball of $ X^{\alpha}$
of radius 
$R(\theta) =  2 KM   \frac{ \Gamma( \frac{1}{2} -\alpha   ) }{\delta^{\frac{1}{2}- \alpha} }  $, for $t \geq T$.


\eproof

\section{Existence of Global Attractors}
\label{attract}

The first step to show the existence of global attractors will be to obtain
a ``contraction property'' of the sets $\Sigma_\theta$, similar to the 
property for rectangles, considered by Smoller  \cite{smoller}.

\begin{lema}
\label{contract}
 Suppose that the hypotheses of Lemma \ref{posinv} hold
and  $\bar\theta\in \R$ satisfy $\bar\theta m_{\epsilon} >\xi$. Then, for any 
$\theta$ there exists a $\bar t$, which can be chosen independently of $\epsilon$, such that 
\[
T_{\epsilon}(t)\Sigma_\theta^{\epsilon}\subset\Sigma_{\bar\theta}^{\epsilon},
\]
for all $t\geq{\bar t}$.
\end{lema}

\medskip

\noindent{\bf Proof:}

\medskip

Let $u\in\Sigma_\theta$. We can suppose without loss of generality that 
$\theta\geq\bar\theta$. Let $\bar v=e^{-t\mu_{\epsilon}}\theta\varphi$, $\underline v=
-\bar v$.
As in Lemma~\ref{posinv}, we can prove that $\bar v$ and $\underline v$ are 
super- and sub-solutions respectively. Thus, using Theorem~\ref{pao} and the 
uniqueness of solution, we have that
\[
\underline v\leq T_{\epsilon} (t)u\leq \bar v,
\]

Therefore  $T_{\epsilon}(t)u$ enters $\Sigma_{\bar\theta}$ after a time depending only on
 $\theta$, and on   the first eigenvalue $\mu_{\epsilon} $ of $A_{\epsilon} $ (and not on the particular solution $u \in \Sigma_{\theta}$). Since $\mu_{\epsilon} $ is bigger than a constant $\mu$, for 
  $\epsilon$ sufficiently small, and 
$\Sigma_{\bar\theta}$
is positively invariant, the result follows.

\eproof

\medskip
%

\begin{teo}
\label{exatr}
 Suppose that the hypotheses of Lemma \ref{posinv} hold. Then
the problem (\ref{abstract_scale}) has a global attractor  $\cala_{\epsilon}$ in $X^{\eta}_{\epsilon}
$. Furthermore 
 $\mathcal{A_{\epsilon}} \subset \Sigma_{\theta}^{\epsilon}$  if $\theta m_{\epsilon}\geq
\xi$.

\end{teo}

\medskip

\noindent{\bf Proof:}

\medskip

Let $V $ be a bounded subset of $ X^{\eta}$, and 
  $\bar\theta\in\R $  be such that $\bar\theta m\geq\xi$. If $u$ is any
element of $X^{\eta}$, it follows from the continuity
of the embedding $X^{\eta}\hookrightarrow C^0(\bar\Omega)$ that
$u\in\Sigma_{\theta}$, 
for some $\theta$ and then, applying Lemma~\ref{contract}, we conclude that 
$T(t)u\in\Sigma_{\bar\theta}$, for $t$ big enough.
From Lemma \ref{bound_theta}, it follows that $V $ enters and remains in a ball
of   $ X^{\alpha} $, with $\alpha > \eta$ of  radius $R(\alpha, \bar{\theta})$,  which does not depend on $V$.   Since this ball is a compact set of
$ X^{\alpha} $, the existence of a global compact  attractor
 $\mathcal{A}$ follows immediately. Furthermore, since $ \Sigma_{\bar{\theta}} $  is positively invariant by
 Lemma
 \ref{posinv} it also follows that $\mathcal{A} \subset \Sigma_{\bar{\theta}} $,
 as claimed.

 \eproof
 \begin{cor} \label{uniform_bound_Linf}
   Suppose that the hypotheses of Lemma \ref{posinv} hold.
 If $\epsilon_0 $ is sufficiently small, the attractor $\cala_{\epsilon}$  is uniformly bounded in
  $ L^{\infty} $, for $0 \leq \epsilon \leq \epsilon_0$. 
\end{cor}

\noindent{\bf Proof:}
 From \eqref{fund_inequality_Lp} and results in \cite{kato},  it follows that the first eigenvalue and eigenfunction of $A_{\epsilon}$ are continuous in  $W^{1,p} $ and, therefore, also in  $L^{\infty}$, Thus  the sets
  $  \Sigma_{{\theta}}^{\epsilon}$ are uniformly bounded in  $ L^{\infty} $ and the result follows from Theorem \ref{exatr}.
\eproof

\section{Uppersemicontinuity of the family of  global attractors}
\label{upper}

  Recall  that a family of subsets  $\mathcal{A}_{\lambda}$ of
 a metric space  $(X,d)$ is
 said to be  {\em upper-semi continuous} at $\lambda = \lambda_0$
 if $\delta(\mathcal{A}_{\lambda}, \mathcal{A}_{\lambda_0}) \to 0$
 as $ \lambda \to \lambda_0$, where
 $\delta(A,B) = \sup_{x\in A} d(x,B) =
  \sup_{x\in A} \inf_{y \in B} d(x,y) $ 
and {\em lower-semicontinuous} if 
 $\delta(\mathcal{A}_{\lambda_0}, \mathcal{A}_{\lambda}) \to 0$
 as $ \lambda \to \lambda_0$.

  To prove the uppersemicontinuity of the family of attractors $A_{\epsilon}$, given by Theorem \ref{exatr}
   in the (fixed) fractional  space 
   $X^{\eta}$, $0 < \eta < \frac{1}{2}$, we will need two main ingredients: the  uniform boundedness of the family and  the continuity of the nonlinear semigroup  $T_{\epsilon}$ with respect to $\epsilon$. This is the content of the next two results. In view of the uniform boundedness of the solutions, proved in Corollary   \ref{uniform_bound_Linf} we may  suppose, without loss of generality, the following hypothesis on the nonlinearites.

   \begin{align} \label{glob_Lip}
    \bullet \  f \textrm{ and } g  & \textrm{ are globally bounded}.  \nonumber  \\
    \bullet \   f \textrm{ and } g & \textrm{ are  globally  Lipschtiz,
      with Lipschitz constants}  L_f  \textrm{ and }  L_g  \ \textrm{respectively},   
  \end{align}
  \begin{lema} \label{uniform_bound_space}
    Suppose that the hypotheses of Lemma \ref{posinv} and \eqref{glob_Lip}  hold.
 If $\epsilon_0 $ is sufficiently small,  the family of attractors  $\cala_{\epsilon}$   given by Theorem \ref{exatr} is uniformly bounded 
   in the (fixed) fractional  space 
  $X^{\eta}$, $0 < \eta < \frac{1}{2}$,  for $0 \leq \epsilon \leq \epsilon_0$. 
\end{lema}

\proof

 Let $b$ be the exponential rate of decay of the  linear semigroup
generated by $A_{\epsilon}$, $ \epsilon$ for $\epsilon $ small, given by 
 Theorem \ref{cont_lin_semigroup}.  Let  $u \in \cala_{\epsilon}$. By the variation of constants formula, Lemma \ref{sector} and Theorem \ref{cont_lin_semigroup}, we obtain  
 \begin{align*}
\| T_{\epsilon}(t)(u)  \|_{\eta} & \leq   
 \|    e^{A_{\epsilon}(t)}  u \|_{\eta}  
 + \int_{0}^{t}  \| e^{A_{\epsilon}(t-s)}
  H_{\epsilon}(T_{\epsilon}(s) u )\|_{\eta} \, ds \\
  & \leq \|    e^{A(t)}  u \|_{\eta}     +  \| \left( e^{A_{\epsilon}(t)} -  e^{A(t)} \right) u \|_{\eta}     
    + \int_{0}^{t}  \| e^{A(t-s)}
  H_{\epsilon}(T_{\epsilon}(s) u )\|_{\eta}  \, d\, s  \\
   &   +  
 \int_{0}^{t}    \|  \left( e^{A_{\epsilon}(t-s)} -  e^{A(t-s)} \right) 
 H_{\epsilon}(T_{\epsilon}(s) u)  \|_{\eta} \, ds   \\
     & \leq  \left( C  e^{-at}  + C(\epsilon)   e^{-bt}   \right)
      \frac{1}{t^{\eta+\frac{1}{2}}} \|u\| 
         + 
  \int_{0}^{t}   C  e^{-a(t-s)} \frac{1}{(t-s)^{\eta}+\frac{1}{2}} 
  \|  H_{\epsilon}(T_{\epsilon}(s) u ) \| \, d \, s \\
   &  +  \int_{0}^{t}    C  e^{ -b (t-s)}  \frac{1}{(t-s)^{\eta + \frac{1}{2}}} \|
  H_{\epsilon}(T_{\epsilon}(s) u ) \| \, ds   \\
   \end{align*}

By \eqref{normH_infty}

 \begin{align*}
\|(H_{\epsilon})_{\eta}(T(s)u_0)\|_{X^{-\frac{1}{2}  }}  &  \leq  M_f   |\Omega|^{\frac{1}{p} } +
 ||\mu||_\infty M_g |{\partial \Omega}|^{\frac{1}{p} }  \,K_1  \\
& \leq 
  \|f\|_{\infty}   |\Omega|^{\frac{1}{p} } +
 ||\mu||_{\infty}  \|g\|_{\infty} |{\partial \Omega}|^{\frac{1}{p} }  \,K_1, 
 \end{align*}
where $K_1$ is a constant of the trace mapping (see Theorem \ref{trace}).
Thus
\begin{align*}
\| T_{\epsilon}(t)(u)  \|_{\eta}    & \leq
  C'   e^{-bt}  \frac{1}{t^{\eta + \frac{1}{2}}}  \|u\|_{\infty}    + C'' \left(  \|f\|_{\infty}   |\Omega|^{\frac{1}{p} } +
 ||\mu||_{\infty}  \|g\|_{\infty} |{\partial \Omega}|^{\frac{1}{p} }  \right)
  \int_{0}^{t}     e^{-b(t-s)} \frac{1}{(t-s)^{\eta}+ \frac{1}{2}}  \, ds,   \\
   \end{align*}
where the constants $C' $ and $C''$ do not depend on $\epsilon$. 

Since  the right hand side is uniformly bounded for $u \in \cala_{\epsilon}, t> 0$ and the attractors are invariant, the result follows immediately.
\eproof

\begin{lema} \label{uniform_cont_space}
  Suppose that the hypotheses of Lemma \ref{uniform_bound_space} hold.
  Then the  map
   \[  (u, \epsilon)\in X^{\eta} \times  [0, \epsilon_0] \mapsto T_{\epsilon} u \in X^{\eta} \] 
 is  continuous at $\epsilon = 0$, uniformly  for $u$ in bounded sets and $ 0 < t  \leq T < \infty$.  
\end{lema}

\proof

Using the variation of constants formula, \eqref{normH_infty}, \eqref{unif_contH_infty} and
 \eqref{LipH}, we obtain

 \begin{align*}
\|  T_{\epsilon}(t)(u) - T(t)(u)  \|_{\eta} & \leq   
 \|    e^{A_{\epsilon}(t)}  u  -  e^{A(t)}  u \|_{\eta} \\
 & + \int_{0}^{t}  \|  \left( e^{A_{\epsilon}(t-s)} - e^{A(t-s)} \right) 
  H_{\epsilon} (T_{\epsilon}(s) u)\|_{\eta} \, ds   \\
   & + 
     \int_{0}^{t}  \|  e^{A(t-s)} \left(  
  H_{\epsilon}(T_{\epsilon}(s) u ) -  H(T_{\epsilon} (s) u )\|_{\eta}    \right)  \| \, ds  \\
  &  +   \int_{0}^{t}  \|  e^{A(t-s)} \left(  
  H(T_{\epsilon}(s) u ) -  H(T (s) u)  \right)\|_{\eta}  \, ds 
       \\       
   & \leq C(\epsilon ) e^{-bt}  \frac{1}{t^{\eta + \frac{1}{2} }}   \|u\|
   +
   \int_{0}^{t}   C (\epsilon)  e^{-b(t-s)} \frac{1}{(t-s)^{\eta}+\frac{1}{2}} 
  \|  H_{\epsilon}(T_{\epsilon}(s) u ) \| \, d \, s \\
   &  +  \int_{0}^{t}    C  e^{ -b (t-s)}  \frac{1}{(t-s)^{\eta + \frac{1}{2}}} \|
  H_{\epsilon}(T_{\epsilon}(s) u ) - H(T_{\epsilon}(s) u)  \| \, ds   \\
  &  +  \int_{0}^{t}    C  e^{ -b (t-s)}  \frac{1}{(t-s)^{\eta + \frac{1}{2}}} \|
  H(T_{\epsilon}(s) u) - H(T(s) u )  \| \, ds   \\
    & \leq C(\epsilon ) e^{-bt} \frac{1}{t^{\eta + \frac{1}{2} }}  \|u\|
   \\
  & +
   \int_{0}^{t}   C (\epsilon)  e^{-b(t-s)} \frac{1}{(t-s)^{\eta}+\frac{1}{2}} 
 \left( \|f\|_{\infty}   |\Omega|^{\frac{1}{p} } +
 ||\mu||_{\infty}  \|g\|_{\infty} |{\partial \Omega}|^{\frac{1}{p} }  \,K_1   \right) \, d \, s \\
   &  +  \int_{0}^{t}    C  e^{ -b (t-s)}  \frac{1}{(t-s)^{\eta + \frac{1}{2}}} \|
 \left(  \|\mu_{\epsilon} - 1 \|_{\infty}  M_g |\partial \Omega|^{\frac{1}{p}  } K_1    \right)  \, ds   \\
   &  +  \int_{0}^{t}    C  e^{ -b (t-s)}  \frac{1}{(t-s)^{\eta + \frac{1}{2}}} \|
  \left( L_g  ||\mu||_\infty  K_1 K_2 +  L_f \right)
    \| T_{\epsilon}(s)  u  - T(s) u \|_{X_{\epsilon}^{\eta}}  \| \, ds   \\
   \end{align*}

   Writting
    \begin{align*}
      A(\epsilon) & : =   C(\epsilon )   \|u\|
   + t^{\eta + \frac{1}{2}}
   \int_{0}^{t}   C (\epsilon)  e^{bs} \frac{1}{(t-s)^{\eta}+\frac{1}{2}} 
 \left( \|f\|_{\infty}   |\Omega|^{\frac{1}{p} }   
 +
 ||\mu||_{\infty}  \|g\|_{\infty} |{\partial \Omega}|^{\frac{1}{p} } \,K_1  \right) \,
     \, d \, s \\
    & +   t^{\eta + \frac{1}{2}}  \int_{0}^{t}    C  e^{ bs}  \frac{1}{(t-s)^{\eta + \frac{1}{2}}} \|
 \left(  \|\mu_{\epsilon} - 1 \|_{\infty}  M_g |\partial \Omega|^{\frac{1}{p}  } K_1    \right)  \, ds  \\
 \\
      B  & : =    C  
  \left( L_g  ||\mu||_\infty  K_1 K_2 +  L_f \right),
 \end{align*}    
  we obtain    

  \begin{align*}
e^{bt} \|  T_{\epsilon}(t)(u) - T(t)(u)  \|_{\eta} & \leq
  A(\epsilon) t^{ - (\eta + \frac{1}{2}) } + 
  B \int_0^t t^{ - (\eta + \frac{1}{2}) }
 e^{ b s}
    \| T_{\epsilon}(s)  u  - T(s) u \|_{X_{\epsilon}^{\eta}}   \, ds   
   \end{align*}

  From the singular Gronwall's inequality, it follows that 
  \begin{align*}
 \|  T_{\epsilon}(t)(u) - T(t)(u)  \|_{\eta} & \leq
  A(\epsilon) M e^{-bt} t^{ - (\eta + \frac{1}{2}) },
   \end{align*}   
   for $0< t  \leq T$, where the constant $M$ depends on $B, \eta $  and $T$, for $u$ in a bounded set
    of $X_{\epsilon}^{\eta}$.
    \eproof

    \begin{teo} \label{upper_semi} Suppose that the hypotheses of Lemma \ref{uniform_bound_space} hold.
      Then the family of
  attractors ${\cala}_{\epsilon}$, given by Theorem \ref{exatr}
  is upper semicontinuous with
respect to $\epsilon$ at $\epsilon= 0$.
\end{teo}
\proof From  Lemma \ref{uniform_bound_space} there exists a bounded set 
$B \subset X^{\eta}$ such that $ \bigcup_{0 \leq \epsilon \leq \epsilon_{0} }  \cala_{\epsilon} \subset B$.  
 Given  $\delta > 0$, there exists  $t_{\delta} >0$  such that  $T(t_{\delta})(B) \subset 
 \cala_0^{\frac{\delta}{2}}$, where
 ${\cala}_{0}^{\frac{\delta}{2}}$ is the
$\frac{\delta}{2}$-neighborhood of ${\cala}_{0}$.

From Lemma \ref{uniform_cont_space},   there exists $\bar{\epsilon}
> 0$ such that 
$| T_{\epsilon}(t_{\delta}) u - T(t_{\delta}) u \|_{X^{\eta}} \leq \frac{\delta}{2} $, for
 every $u \in B$ and $ 0 \leq \epsilon \leq \bar{\epsilon}$. It follows that 
$ T_{\epsilon}(t_{\delta}) B \subset \cala_0^{\delta}$. In particular,   
$ T_{\epsilon}(t_{\delta}) \cala_{\epsilon} \subset \cala_0^{\delta}$.
Since $\cala_{\epsilon}$ is invariant under $T_{\epsilon}$, we conclude that 
$  \cala_{\epsilon} \subset \cala_0^{\delta}$, for $0\leq \epsilon \leq \bar{\epsilon}$, thus proving the claim.
\eproof

From the semicontinuity of attractors, we can easily prove the corresponding
 property for the equilibria.
\begin{cor}\label{upperequil}
 Suppose the hypotheses of Theorem \ref{upper_semi} hold.
   Then the family of sets of equilibria
   $ {\{} E_{\epsilon } \,   {|} \, 0\leq \epsilon \leq \epsilon_0  {\}}$, of
     the problem  (\ref{abstract_scale}) is uppersemicontinuous in $X^\eta$.
\end{cor}
\proof
The result is well-known, but we sketch a proof here for completeness.
Suppose $u_{n} \in \cala_n$, with $\displaystyle\lim_{n\to \infty}\epsilon_n= 0$. We choose an arbitrary subsequence and still call it   $(u_{n})$, for simplicity. It is enough to show that, there exists a subsequence   $(u_{n_k})$, which converges to a point $u_0 \in E_0$. Since $(u_{n}) \to \mathcal{A}_{0}$,  there exists
$(v_n) \in \cala_0$ with $ \|u_n - v_n\|_{\eta} \to 0 $. Since    $\mathcal{A} {_0}$ is compact,
there
exists  a subsequence   $(v_{n_k})$, which converges to a point $u_0 \in
\mathcal{A}_0$, so also  $(u_{n_k}) \to \mathcal{A}_0 $.
Now, since the flow $T_\epsilon(t)$ is continuous in $\epsilon$
we have, for
any $t>0$
\[ u_{n_k} \to u_0 \Leftrightarrow  T_{\epsilon_{n_k}}(t) u_{n_k} \to
  T_{0}(t) u_0 \Leftrightarrow   u_{n_k} \to
  T_{0}(t) u_0. \]
  Thus, by uniqueness of the limit, $T_{0}(t) u_0 =u_0$, for any $t> 0$, so
 $u {_0} \in E_0. $ 
\eproof

\section{Lowersemicontinuity}

\label{lower}
For the lower semicontinuity  we will need to 
assume the following  additional properties for the  nonlinearities.
\vspace{3mm}

\begin{align} \label{boundfg}
 f  \textrm{ and }   g  \textrm{ are in }    C^1(\R,\R)  \textrm{ with bounded derivatives }.  
\end{align}



\begin{lema}
 \label{FGateaux}
 Suppose that $\eta$ and $p$  are such that \eqref{hip_inclusion} holds and $f$ satisfies  (\ref{boundfg}). 
 Then 
the operator 
 $F :X^{\eta}\times  \R  {\rightarrow} X^{-\frac{1}{2}}$ given by 
 (\ref{Fh}) is Gateaux differentiable with respect to
 $u$, with Gateaux differential
 $\displaystyle{\frac{\partial F}{\partial u}(u,\epsilon)w}$ given by 
 \begin{equation}\label{FGateaux_form}
\left\langle \frac{\partial F}{\partial u}(u,\epsilon)w\,,\,
\Phi\right\rangle_{-\frac{1}{2},\frac{1}{2}} = \displaystyle\int_\Omega f^{\,'}(u)w\,\Phi\,dx\,,
\end{equation}
for all  $w \in X^\eta$ and $\Phi \in X^{\frac{1}{2}}$.
\end{lema}
\proof
 Observe first that $F(u, {\epsilon})$ is well-defined, since the conditions of Lemma \ref{Flip} are met.

It is clear that $\displaystyle\frac{\partial F}{\partial u}(u,\epsilon)$ is linear. We now show that it is bounded. In fact we have,
 for all 
 $u,w \in X^{\eta} $  and
 $\Phi \in X^{-\frac{1}{2} }  = W^{1,q}$
\begin{eqnarray*}
\left|\left\langle \frac{\partial F}{\partial u}(u,\epsilon)w\,,\,
\Phi\right\rangle_{-\frac{1}{2},\frac{1}{2}} \right|  &\leq & \displaystyle\int_\Omega 
|\,f^{\,'}(u)\,|\, |\,w\,| \,|\,\Phi\,|\, dx\ \\
&\leq &  \|\,f'\,\|_{\infty}\displaystyle\int_\Omega |\,w\,| \, |\,\Phi\,|\,dx \\
&\leq &  \|\,f'\,\|_{\infty}\displaystyle \|\,w\,\|_{L^p(\Omega)} \,
 \|\,\Phi\,\|_{L^q(\Omega)}\,dx \\
&\leq &  \|\,f'\,\|_{\infty}\displaystyle \|w\|_{ X^{\eta}} \,
 \|\,\Phi\,\|_{ X^{\frac{1}{2}}}\,dx\,, 
\end{eqnarray*}
where $\|f'\|_{\infty}= \sup \{f'(x)  \,|\, x \in \R \} $.  This proves boundedness. 

Now, we  have,  for all 
 $u,w \in X^{\eta} $  and
 $\Phi \in X^{\frac{1}{2}}$
\begin{eqnarray*}
&& \left|\frac{1}{t}\left\langle F(u + tw,\epsilon) - F(u,\epsilon) - t  \frac{\partial F}{\partial u}(u,\epsilon) w,\Phi\right\rangle_{-\frac{1}{2},\frac{1}{2}}\right| \\
&\leq& \frac{1}{|t|}\displaystyle \int_\Omega \big|\,[\,f(u + tw) - f(u)
  - tf^{\,'}(u)w\,]\,\Phi\,\big|\,dx \\
&\leq &  \frac{1}{|t|} \left(\displaystyle\int_\Omega \big|f(u + tw) - f(u) - tf^{\,'}(u) w\big|^{p}dx\right)^\frac{1}{p}||\Phi||_{X^{\frac{1}{2}}} \\
& \leq & 
 \left(\,\displaystyle \underbrace{\int_\Omega \big|\,
 \left(f'(u + \bar{t}w ) -  f^{\,'}(u)\right) w\,\big|^{\,p} \,dx\,
}_{(I)}\right)^\frac{1}{p}||\,\Phi\,||_{X^{\frac{1}{2}}}, 
\end{eqnarray*}
where $  0 \leq \bar{t} \leq t$.  
Since $f'$ is bounded and continuous, the integrand of $(I)$
  is bounded by an integrable function and goes to $0$ as $t \to 0$.
 Thus, the integral $(I)$ goes to $ 0$  as $t \to 0$, from Lebesgue's Dominated
Convergence Theorem. It follows that \newline
$ \displaystyle{\lim_{t \to 0} \frac{ F(u + tw {,\epsilon}) - F(u {,\epsilon})}{t}  = \frac{\partial F}{\partial u}(u,\epsilon) w \ \textrm{ in}  \  X^{-\frac{1}{2}},}$
 for all  $u,w \in X^{\eta} $; so $F$ is Gateaux differentiable with 
 Gateaux differential given by (\ref{FGateaux_form}).
\eproof

We now want to prove that the Gateaux differential of $F(u, \epsilon)$
is continuous in $u$. Let us denote by
 $\mathcal{B}(X, Y)$  
the space of
  linear bounded operators from $X$ to $Y$.
 We will need the following result, whose simple proof is omitted.

 \begin{lema}\label{strong_uniform_operators} 
Suppose  $X,Y$ are Banach spaces and $ T_n : X \to Y$ is a sequence
 of linear operators converging strongly to the linear operator
 $T:X \to Y$. Suppose also that 
 $X_1 \subset X$  is a Banach space, the inclusion 
 $i: X_1 \hookrightarrow X$
 is compact and let $ \widetilde{T}_n = T_n \circ i$ 
 and  $ \widetilde{T} = T \circ i$.   Then
  $\widetilde{T}_n \to  \widetilde{T} $ uniformly for $x$ in a 
 bounded subset of $X_1$ (that is, in the or norm  of 
 $\mathcal{B}(X_1,Y )$).
\end{lema}

\begin{lema}\label{FGateaux_cont}
 Suppose that $\eta$ and $p$  are such that \eqref{hip_inclusion} holds and $f$ satisfies  (\ref{boundfg}). 
Then 
the Gateaux differential of $F(u,\epsilon)$, with respect to $u$ is 
 continuous in $u$, that is, the map 
  $ u \mapsto \displaystyle\frac{\partial F}{\partial u}( u,\epsilon)
 \in \mathcal{B}(X {^\eta}, X^{-\frac{1}{2}})$
 is continuous. 
\end{lema}
\proof Let  $u_n$ be a sequence converging to
 $u$ em $X^{\eta}$, and choose $0<\widetilde{\eta}< \eta$, such that the hypotheses still hold.
 Then, we have for any
 $\Phi \in X^{\frac{1}{2}}$ and  $w \in X^{\widetilde{\eta}}$:
\begin{eqnarray*}  
\bigg|\left\langle \left(  \frac{\partial F}{\partial u}(u_n, {\epsilon}) -  \frac{\partial F}{\partial u}(u, {\epsilon})\right)w\,,\,\Phi \right\rangle_{-\frac{1}{2}\,,\, \frac{1}{2}}\bigg|
&\leq& \displaystyle\int_\Omega \bigg|\,\big(\,f^{\,'}(u) - f^{\,'}(u_n)\,\big)w\,\Phi\,\bigg|\,dx \nonumber \\
&\leq& \bigg(\displaystyle\int_\Omega \big|\big(f^{\,'}(u) - f^{\,'}(u_n)\big)w\big|^{p}dx\bigg)^\frac{1}{p}\bigg(\displaystyle\int_\Omega |\Phi\big|^{q}dx\bigg)^\frac{1}{q} \nonumber \\
&\leq& \,\bigg(\displaystyle\underbrace{\int_\Omega \big|\big(f^{\,'}(u) - f^{\,'}(u_n)\big)w\big|^{p}dx}_{(I)}\bigg)^\frac{1}{p}\,||\Phi||_{ {X^{\frac{1}{2}}}}\,, 
\end{eqnarray*}

\par Now, the integrand in $(I)$ is bounded by the
 integrable function $\,||\,f^{\,'}||_\infty^{\,p}\,w^{\,p}$ and
 goes to $0$ a.e. as   $u_n \to u$ in $X^\eta$.
 Therefore the sequence of operators
$\displaystyle\frac{\partial F}{\partial u}( u_n,\epsilon)$ converges strongly
in the space $\mathcal{B}(X^{\widetilde{\eta}}, X^{-\frac{1}{2}})$  to the operator  $ \displaystyle\frac{\partial F}{\partial u}( u,\epsilon)$.
 From Lemma \ref{strong_uniform_operators} 
 the convergence holds in the norm of  $\mathcal{B}(X^{\eta}, X^{-\frac{1}{2}})$,
 since  $X^{\eta}$ is compactly embedded in $X^{\widetilde{\eta}}$.
\eproof

\begin{lema}
 \label{GGateaux}
  Suppose that $\eta$ and $p$  are such that \eqref{hip_inclusion} holds and $g$ satisfies  (\ref{boundfg}).
 Then
the operator 
 $G :X^{\eta}\times  \R  {\rightarrow} X^{-\frac{1}{2}}$ given by 
 (\ref{Gh}) is Gateaux differentiable with respect to
$u$, with Gateaux differential
\begin{equation}\label{GGateaux_form}
\left\langle \frac{\partial G}{\partial u}(u, {\epsilon})w\,,\,\Phi\right\rangle_{ {-\frac{1}{2}, - -\frac{1}{2}}} = \displaystyle\int_{\partial\Omega} g^{\,'}(\gamma(u))\gamma(w)\,\gamma(\Phi)\,\left|\displaystyle\frac{J_{\partial\Omega}h_\epsilon}{Jh_\epsilon}\right|\,d\sigma(x)\,,
\end{equation}
for all  $w \in X^\eta$ and $\Phi \in X^{\frac{1}{2}}$.
\end{lema}
\proof
Observe first that $G(u, {\epsilon})$ is well-defined, since the conditions of Lemma
\ref{Gbem} are met.

It is clear that $\displaystyle\frac{\partial G}{\partial u}(u,\epsilon)$ is linear. We now show that it is bounded. In fact we have,
 for all 
 $u,w \in X^{\eta} $  and
 $\Phi \in X^{\frac{1}{2}}$
 \begin{eqnarray*}
   \left| \left\langle \frac{\partial G}{\partial u}(u, {\epsilon})w\,,\,\Phi\right\rangle_{-\frac{1}{2}, \frac{1}{2}} \right| &  = &  \left| \displaystyle\int_{\partial\Omega} g^{\,'}(\gamma(u))\gamma(w)\,\gamma(\Phi)\,\left|\displaystyle\frac{J_{\partial\Omega}h_\epsilon}{Jh_\epsilon}\right|\,d\sigma(x)\, \right| \\
   &\leq &
   \|\mu \|_{\infty} \, \|g'\|_{\infty}
   \displaystyle\int_{\partial\Omega} | \gamma(w)|\, |\gamma(\Phi) |\,
   \,d\sigma(x)\,
   \\
   &\leq &  \|\mu \|_{\infty} \, \|g'\|_{\infty}
   \displaystyle \|\gamma(w)\|_{L^p(\partial \Omega)} \,
 \|\,\gamma(\Phi)\,\|_{L^q(\partial \Omega)}\\
 &\leq &
  K_1 K_2  \|\mu \|_{\infty} \, \|g'\|_{\infty}
   \displaystyle \|w\|_{X^{\eta}}} \,
 \|\,\Phi\,\|_{X^{\frac{1}{2}} \,,
\end{eqnarray*}
 where $\|g'\|_{\infty}= \sup \left\{g'(x)  \,|\, x \in \R \right\} $,
$\|\mu \|_{\infty} =   \sup \left\{ |\mu(x, \epsilon)| \, | \, 
 x\in \partial \Omega \right\} =  \sup \left\{
 \displaystyle\frac{J_{\partial\Omega}h_\epsilon}{Jh_\epsilon} (x) \, \bigg| \, 
 x\in \partial \Omega \right\}$
 and $K_1$, $K_2$ are embedding constants given by Theorem \ref{trace}. This proves boundedness. 

Now, we  have,  for all 
 $u,w \in X^{\eta} $  and
 $\Phi \in X^{\frac{1}{2}}$
\begin{eqnarray*}
 & & \left| \frac{1}{t} \left\langle G(u + tw {,\epsilon}) - G(u {,\epsilon}) - t  \frac{\partial G}{\partial u}(u,\epsilon) w\,,\,\Phi\right\rangle_{-\frac{1}{2},\frac{1}{2}}\right|  \\ 
 &\leq&  \frac{1}{|t|}\displaystyle \int_{\partial  {\Omega}} \left|\,\left[\,g(\gamma(u + tw)) -g(\gamma(u))
 - tg'(\gamma(u))\right]\gamma(w)\,\right|\,  {\left|\gamma(\Phi)\right|}\,\left| \frac{J_{\partial\Omega}h_\epsilon}{Jh_\epsilon} \right|\, {d\sigma(x)}   \\
  & \leq& K_1  \|\mu \|_{\infty} 
 \frac{1}{|t|}\displaystyle \left\{ \int_{\partial  {\Omega}} \left|\,\left[\,g(\gamma(u + tw)) -
 g(\gamma(u))
  - tg'(\gamma(u))\,\right]\gamma(w)\right|^p \, {d\sigma(x)} \right\}^{\frac{1}{p}} 
  \| \Phi \|_{X^{\frac{1}{2}}}    \\  
 &\leq& K_1 \|\mu \|_{\infty}  \displaystyle \left\{ \underbrace{\int_{\partial  {\Omega}}
 \left|\,\left[\,g'(\gamma(u + \bar{t}w)) 
  - g'(\gamma(u))\,\right]\gamma(w)\,\right|^p \, {d\sigma(x)}}_{(I)} \right\}^{\frac{1}{2}} 
  \| \Phi \|_{X^{\frac{1}{2}}} \,,  
\end{eqnarray*}
where  $K_1$  is
the  embedding constant given by Theorem \ref{trace} and Lemma \ref{inclusion} and $  0 \leq \bar{t} \leq t$.  
Since $g'$ is bounded and continuous, the integrand of $(I)$
  is bounded by an integrable function and goes to $0$ as $t \to 0$.
 Thus, the integral $(I)$ goes to $ 0$  as $t \to 0$, from Lebesgue's Dominated
Convergence Theorem. It follows that 
$ \displaystyle{\lim_{t \to 0} \frac{ G(u + tw {,\epsilon}) - G(u {,\epsilon})}{t}  =
 \frac{\partial G}{\partial u}(u,\epsilon) w \ \textrm{ in}  \  X^{-\frac{1}{2}},}$
 for all  $u,w \in X^{\eta} $; so $G$  is Gateaux differentiable with 
 Gateaux differential given by (\ref{GGateaux_form}).
\eproof

\begin{lema}\label{GGateaux_cont}
 Suppose that $\eta$ and $p$  are such that \eqref{hip_inclusion} holds  and $g$ satisfies  (\ref{boundfg}).
Then
the Gateaux differential of $G(u,\epsilon)$, with respect to $u$ is 
 continuous in $u$ (that is, the map 
  $ u \mapsto \displaystyle\frac{\partial G}{\partial u}( u,\epsilon)
 \in \mathcal{B}(X^{\eta}, X^{-\frac{1}{2}})$
 is continuous) and uniformly continuous in $\epsilon$ for $u$ in bounded sets
 of $X^{\eta}$ and
 $0\leq \epsilon \leq \epsilon_0 <1$. 
\end{lema}

\proof Let $0\leq \epsilon \leq \epsilon_0$,  $u_n$ be a sequence converging to
 $u$ em $X^{\eta}$, and choose $ 0< \widetilde{\eta} < \eta$, still satisfying the hypotheses. 
 Then, we have for any
 $\Phi \in X {^{\frac{1}{2}}}$ and  $w \in X^{\widetilde{\eta}}$:
\begin{eqnarray*}
& &\bigg|\,\left\langle \,\left(\,  \frac{\partial G}{\partial u}(u_n,\epsilon) -  \frac{\partial G}{\partial u}(u,\epsilon)\,\right)w\,,\,\Phi \,\right\rangle_{-\frac{1}{2}\,,\, \frac{1}{2}}\,\bigg| \\
& \leq & \displaystyle\int_{\partial\Omega}
\left| \left(g'(\gamma(u)) - g'(\gamma(u_n))\right)
\gamma(w)\,\gamma(\Phi)\right|\,\left|\displaystyle\frac{J_{\partial\Omega}h_\epsilon}{Jh_\epsilon}\right|\,d\sigma(x)
\nonumber \\
& \leq & \|\mu_{\epsilon} \|_{\infty} \left\{ \displaystyle\int_{\partial\Omega}
\left| (g'(\gamma(u)) - g'(\gamma(u_n))
\gamma(w)\,\,\right|^p\,d\sigma(x)\, \right\}^{\frac{1}{p}}
  \left\{ \displaystyle\int_{\partial\Omega}
\left| \,\gamma(\Phi)\,\right|^q\,d\sigma(x)\, \right\}^{\frac{1}{q}} 
\nonumber \\
 &\leq &K_1 \|\mu_{\epsilon}\|_{\infty} \left\{ \displaystyle \underbrace{\int_{\partial\Omega}
\left| (g'(\gamma(u)) - g'(\gamma(u_n))
\gamma(w)\,\,\right|^p\,d\sigma(x)}_{(I)}\, \right\}^{\frac{1}{p}}
 \|\Phi\|_{ {X^{\frac{1}{2}} }}\,,
\nonumber 
\end{eqnarray*}
where $K_1$ is the  constant due to continuity of the trace map from
$X^{\frac{1}{2}}$ into
 $L^2(\partial \Omega)$, as in Lemma \ref{Gbem}.

\par Now, the integrand in $(I)$ is bounded by the
 integrable function $||\,g'\,||_\infty^2 \left| 
\gamma(w)\right|^2$ and
goes to $0$ a.e. as   $u_n \to u$  {in} $X^\eta$.
 Therefore the sequence of operators
 $ \displaystyle\frac{\partial G}{\partial u}( u_n,\epsilon)$ converges strongly
in the space $\mathcal{B}(X^{\widetilde{\eta}}, X^{-\frac{1}{2}})$  to the operator  $ \displaystyle\frac{\partial G}{\partial u}( u,\epsilon)$.
 From Lemma \ref{strong_uniform_operators} 
 the convergence holds in the norm of  $\mathcal{B}(X^{\eta}, X^{-\frac{1}{2}})$,
 since  $X^{\eta}$ is compactly embedded in $X^{\widetilde{\eta}}$ (see \cite{He2}).

 Finally, if $0\leq \epsilon_1 \leq \epsilon_2 <\epsilon_0$, we have
 for any
 $\Phi \in  {X^{\frac{1}{2}}}$ and  $w \in X^{\eta}$:
\begin{eqnarray*} 
& &\bigg|\,\left\langle \,\left(\,  \frac{\partial G}{\partial u}(u,\epsilon_1) -  \frac{\partial G}{\partial u}(u,\epsilon_2)\,\right)w\,,\,\Phi \,\right\rangle_{-\frac{1}{2}\,,\, \frac{1}{2}}\,\bigg| \\
&  \leq & \displaystyle\int_{\partial\Omega}
\left|\,g'(\gamma(u))
\gamma(w)\,\gamma(\Phi)\,\right|\,\left|\displaystyle \mu_{\epsilon_1} -
\mu_{\epsilon_2}\right|\,d\sigma(x)\,,
\nonumber \\
& \leq&  \|\mu_{\epsilon_1} -
\mu_{\epsilon_2}  \|_{\infty} \left\{ \displaystyle\int_{\partial\Omega}
\left| g'(\gamma(u))
\gamma(w)\,\,\right|^p\,d\sigma(x)\, \right\}^{\frac{1}{p}}
  \left\{ \displaystyle\int_{\partial\Omega}
\left| \,\gamma(\Phi)\,\right|^q \,d\sigma(x)\, \right\}^{\frac{1}{q}} 
\nonumber \\
 & \leq & K_1 K_2\|\,g'\,\|_{\infty}\|\|\,w\, \|_{ {X^\eta}}
 \|\Phi\|_{ {X^{\frac{1}{2}}}}
  \|\mu_{\epsilon_1}-\mu_{\epsilon_2}\|_{\infty},
 \nonumber
\end{eqnarray*}
where $K_2$ is the  constant due to continuity of the trace map from
$X^{\eta}$ into
 $L^q(\partial \Omega)$, as before.
 This proves  uniform continuity in $\epsilon$.
\eproof

 \begin{lema}\label{Hfrechet}
  
    Suppose that $\eta$ and $p$  are such that \eqref{hip_inclusion} holds and $f$ and  $g$ satisfy (\ref{boundfg}).
 Then, 
  the map
 $ {(}H {_\epsilon)_-\frac{1}{2}}=  {(}F {_\epsilon)_-\frac{1}{2}} +  {(}G {_\epsilon)_-\frac{1}{2}}  :X^{\eta}\times  \R \mapsto X^{-\frac{1}{2}}$ given by 
 (\ref{defH}) is continuously Fr\'echet differentiable with respect to
  $u$ and the derivative $\displaystyle\frac{\partial G}{\partial u}$ is uniformly continuous with respect to $\epsilon$, for $u$ in bounded sets
  of $X^{\eta}$ and $0\leq \epsilon \leq\epsilon_0 < 1$.
\end{lema}
\proof
The proof follows from Lemmas \ref{FGateaux_cont}, \ref{GGateaux_cont}
and Proposition 2.8 in \cite{Rall}.
\eproof 

We now prove lower semicontinuity for the equilibria.

\begin{teo}\label{equicon}
   {If $f$} and $g$   satisfy the conditions of Theorem  \ref{posinv}
  and
  also (\ref{boundfg}), the equilibria
  of  (\ref{abstract_scale}) with $\epsilon = 0$ are all hyperbolic  
  and $\frac{1}{4}<\eta< \frac{1}{2}$, then the  family of sets of equilibria
  $\{ E_{\epsilon} \, | \, 0 \leq \epsilon <\epsilon_0 \}$ of
  (\ref{abstract_scale}) is
  lower semicontinuous in $X^{\eta}$ at $\epsilon = 0$.
\end{teo}
\proof A point $e \in X^{\eta} $ is an equilibrium of (\ref{abstract_scale})
if and only if it is a root of the map
$$
\begin{array}{rlc}
Z: W^{1,p}(\Omega) \times  {\R}& \longrightarrow &X^{-\frac{1}{2}} \, \\
(u\,,\,\epsilon)& \longmapsto & (A_{ {\epsilon}})_{-\frac{1}{2}}(u) + (H_{\epsilon})_{-\frac{1}{2}}(u)\,,
\end{array}
$$

By Lemma \ref{Hfrechet} the map $ {(}H_\epsilon {)_{-\frac{1}{2}}}: X^{\eta} \to X^{-\frac{1}{2}}$ is continuously Fr\'echet differentiable with
respect to $u$ and by Lemmas \ref{Glip} and \ref{Flip} it is also continuous in $\epsilon$
if $\eta= \frac{1}2 - \delta$, with $\delta>0$ is sufficiently small.
Therefore, the same holds if $\eta = \frac{1}{2}$.

The map  $A_\epsilon= -h_\epsilon^{*} \Delta_{\Omega_\epsilon} h_\epsilon^{*} \,+\, aI$
is a bounded  linear operator from  $W^{1,p}(\Omega)$ to $X^{-\frac{1}{2}}$.
It is also
 continuous in $\epsilon$ since it is analytic as a function of
$h {_\epsilon} \in Diff {^1}(\Omega)$ and
 $ h_\epsilon$ is continuous in  $\epsilon$.
 
 Thus, the map $Z$ is continuously differentiable in $u$ and continuous in
  $\epsilon$. 
 The derivative of $\displaystyle\frac{\partial Z}{\partial u}(e, 0)$
 is an isomorphism by hypotheses.
 Therefore,  the Implicit Function Theorem apply, implying that the
 zeroes of $Z(\cdot, \epsilon)$ are given by a continuous function
 $ e(\epsilon)$. This proves the claim. \eproof

 To prove the lower semi continuity of the attractors, we also need  the continuity of local unstable manifolds at equilibria.

\begin{teo}\label{manifcont}
Suppose that $\eta$ and $p$  are such that \eqref{hip_inclusion} holds  and $f$ and  $g$ satisfy (\ref{boundfg}).,
   $u_0$ is an equilibrium of   (\ref{abstract_scale}) with $\epsilon = 0$,
   and for each $\epsilon>0$ sufficiently small, let $u {_\epsilon}$
   be the unique equilibrium of  (\ref{abstract_scale}), whose existence
   is asserted by Corollary \ref{upperequil} and Theorem \ref{equicon}.
   Then, for $\epsilon$ and $\delta$ sufficiently small, there exists a
   local unstable manifold
$
W_{\rm loc}^u(u_{\epsilon}) 
$ of $u_{\epsilon}$, and if we denote
$ W_{\delta}^u(u_{\epsilon}) =\{ w \in  W_{\rm loc}^u(u_{\epsilon})  \ | \
\|w-u_{\epsilon} \|_{X^{\eta}} < \delta  \}, then$
\[
  -\frac{1}{2} \Big(W_{\delta}^u(u_{\epsilon}),W_{\delta}^u( u_0) \Big) \quad \textrm{and} \quad
   -\frac{1}{2} \Big(W_{\delta}^u(u_{0}),W_{\delta}^u( u_{\epsilon}) \Big)
\]
approach zero as $\epsilon \to 0$, where
  $-\frac{1}{2}(O,Q)=\displaystyle\sup_{o \in O} \inf_{q \in Q}
\|q-o\|_{ {X^{\eta}}}$ for $O$, $Q\subset X^{\eta}$.
\end{teo}
 \proof
 Let $H_{\epsilon}(u)=H(u,\epsilon)$ be the map defined by (\ref{defH}) and $u_{\epsilon}$ a hyperbolic equilibrium of (\ref{abstract_scale}). Since $H(u,\epsilon)$ is differentiable by Lemma \ref{Hfrechet},
 it follows that  $H_{\epsilon}(u_{\epsilon}+w , \epsilon)=
 H_{\epsilon}(u_{\epsilon},\epsilon)
   + H_u(u_{\epsilon} , \epsilon)w + r(w,\epsilon)= A_{\epsilon}u_{\epsilon} + H_u(u_{\epsilon} , \epsilon)w + r(w , \epsilon)$,
 with $r(w,\epsilon)=o(\|w\|_{X^\eta})$, as $\|w\|_{X^\eta} \to 0$.
 The claimed result was proved in \cite{PP}, assuming the following properties of $H_{\epsilon}$:
 
 \begin{itemize}
\item[a)]  $||\,r(w,0)-r(w,\epsilon)\,||_{X^{-\frac{1}{2}}} \leq C({\epsilon})$,
  with  $C({\epsilon}) \to 0 \textrm{ when } \epsilon \to 0$, uniformly for $w$ in a neighborhood of $0$ in $X^{\eta}$.
\item[b)] $||\,r(w_1,\epsilon)-r(w_2,\epsilon)\,||_{X^{-\frac{1}{2}}} \leq k(\rho) ||\,w_1-w_2\,||_{X^{\eta}}$ , for $||\,w_1\,||_{X^{\eta}}\leq \rho$, $||\,w_2\,||_{X^{\eta}}\leq \rho$, with $k(\rho) \to 0$ when $\rho \to  0^+$  and $k(*)$ is non decreasing.
\end{itemize}

 Property a) follows from easily from the fact that both $H(u,\epsilon)$ and
 $H_u(u,\epsilon)$ are
 uniformly continuous in $\epsilon$ for $u$ in bounded sets
 of $X^{\eta}$, by Lemmas \ref{Glip},  {\ref{Flip}}  and  \ref{Hfrechet}.
 It remains to prove property b).
 
 If  $w_1,w_2 \in X^{\eta}$ and  $\epsilon \in [0, \epsilon_0]$,  with
  $0 < \epsilon_0 <1 $  small enough, we have 
\begin{eqnarray}
||\,r(w_1\,,\,\epsilon)-r(w_2\,,\,\epsilon)\,||_{X^{-\frac{1}{2}}} &=& ||\,H(u_{\epsilon} + w_1\,,\,\epsilon) - H(u_{\epsilon}\,,\,\epsilon) - H_u(u_{\epsilon}\,,\,\epsilon)w_1 \nonumber\\ &&-\,
H(u_{\epsilon} + w_2\,,\,\epsilon) + H_{\epsilon}(u_{\epsilon}\,,\,\epsilon) + H_u(u_{\epsilon}\,,\,\epsilon)w_2\,||_{X^{-\frac{1}{2}}} 
\nonumber\\ 
& \leq& ||\,F(u_{\epsilon} + w_1\,,\,\epsilon) - F(u_{\epsilon}\,,\,\epsilon) - F_u(u_{\epsilon}\,,\,\epsilon)w_1 \label{7}\\ 
&&-\, F(u_{\epsilon} + w_2\,,\,\epsilon) + F(u_{\epsilon}\,,\,\epsilon) + F_u(u_{\epsilon}\,,\,\epsilon)w_2\,||_{X^{-\frac{1}{2}}}\nonumber \\
&&+\, ||\,G(u_{\epsilon} + w_1\,,\,\epsilon) -G(u_{\epsilon}\,,\,\epsilon) - G_u(u_{\epsilon}\,,\,\epsilon)w_1 \label{8}\\ 
&&-\, G(u_{\epsilon} + w_2\,,\,\epsilon) + G(u_{\epsilon}\,,\,\epsilon) + G_u(u_{\epsilon}\,,\,\epsilon)w_2\,||_{X^{-\frac{1}{2}}} \,.\nonumber
\end{eqnarray}
 We first estimate (\ref{7}). Since $f'$ is  bounded by (\ref{boundfg}), we have

\begin{align*}
& \bigg|\,\left\langle \,F(u_{\epsilon} + w_1\,,\,\epsilon) - F(u_{\epsilon}\,,\,\epsilon) - F_u(u_{\epsilon}\,,\,\epsilon)w_1 \right.
  \left. -F(u_{\epsilon} + w_2\,,\,\epsilon) + F(u_{\epsilon}\,,\,\epsilon) + F_u(u_{\epsilon}\,,\,\epsilon)w_2 \,,\,\Phi\,\right\rangle_{ {-\frac{1}{2}, \frac{1}{2}}}
  \,\bigg| \\
&\leq
\displaystyle\int_{\Omega} \left|\,[\,f(u_{\epsilon}+w_1)-f(u_{\epsilon})- f'(u_{\epsilon})w_1
  -f(u_{\epsilon}+w_2)+ f(u_{\epsilon})+ f'(u_{\epsilon})w_2\,]\,\Phi\, \right|\,dx\,
\nonumber\\
& =
\displaystyle\int_{\Omega} \left|\,[\,f'(u_{\epsilon}+ \xi_x)-f'(u_{\epsilon})\,](w_1(x)-w_2(x))\,\Phi \,\right|\,dx \nonumber\\
& \leq
K_1  \displaystyle \left\{\int_{\Omega} \left|\,[\,f'(u_{\epsilon}+ \xi_x)-f'(u_{\epsilon})\,]^p
(w_1(x)-w_2(x))^p\,\right|\,dx\right\}^{\frac{1}{p}} \|\Phi \|_{X^{\frac{1}{2}} }\,,
\nonumber\\
& \leq  K_1  K_2 \displaystyle \left\{\int_{\Omega} \left|\,[\,f'(u_{\epsilon}+ \xi_x)-f'(u_{\epsilon})\,]^p
\,\right|\,dx\right\}^{\frac{1}{p}}   \| w_1 -w_2 \|_{X^{\eta}   }  \cdot \|\Phi \|_{X^{\frac{1}{2}} }\,,
\end{align*}
where   $K_1$ is  the embedding constant 
of   $X^{\frac{1}{2}}$  {into} $L^q(\Omega)$, $K_2$ is the embedding constant of 
 $  X^{\eta} $ in $L^{\infty} (\Omega) $ and 
$ w_1(x) \leq \xi_x \leq w_2(x)$ or $ w_2(x) \leq \xi_x \leq w_1(x)$. 
Therefore, we have
$$
\begin{array}{lll}
&&||\,F(u_{\epsilon} + w_1,\epsilon) - F(u_{\epsilon},\epsilon) - F_u(u_{\epsilon},\epsilon)w_1 -
F(u_{\epsilon} + w_2,\epsilon) + F(u_{\epsilon},\epsilon) + F_u(u_{\epsilon},\epsilon)w_2||_{X^{-\frac{1}{2}}} \\
&\leq& K_1K_2\left\{ \displaystyle\int_{\Omega} [\,f'(u_{\epsilon}+ \xi_x)-f'(u_{\epsilon})\,]^{\,p}dx \right\}^\frac{1}{p}
||\,w_1-w_2\,||_{X^{\eta}}\,.
\end{array}
$$
Now the integrand above is bounded by  $2^{\,p}||f'||_\infty^{\,2p}$
 and goes a.e. to $0$ as $\rho \to 0$, since   $||\,w_1\,||_{X^{\eta}}\leq \rho$,  $||\,w_2\,||_{X^{\eta}}\leq \rho$  and  $w_1(x)\leq \xi_x\leq w_2(x)$. 
Thus, the integral goes to $0$ by Lebesgue's bounded convergence Theorem.
 
\quad We now estimate (\ref{8}):

\begin{align*}
& \bigg|\,\left\langle G(u_{\epsilon} + w_1\,,\,\epsilon) - G(u_{\epsilon}\,,\,\epsilon) - G_u(u_{\epsilon}\,,\,\epsilon)w_1 
   -G(u_{\epsilon} + w_2\,,\,\epsilon) + G(u_{\epsilon}\,,\,\epsilon) + G_u(u_{\epsilon}\,,\,\epsilon)w_2 \,,\,\Phi\,\right\rangle_{ {-\frac{1}{2}, \frac{1}{2}}}\,\bigg| \\
& \leq
\displaystyle\int_{\partial\Omega} \bigg|\,[\,g(\gamma(u_{\epsilon}+w_1)) - g(\gamma(u_{\epsilon}))- g'(\gamma(u_{\epsilon}))w_1 \\
&    \quad \quad \quad \quad \quad  \quad \quad -g(\gamma(u_{\epsilon}+w_2)) + g(\gamma(u_{\epsilon}))+ g'(\gamma(u_{\epsilon}))w_2\,]
\,\gamma(\Phi)\gamma\left(\,\left|\,\displaystyle\frac{J_{\partial\Omega} {h_\epsilon}}{J {h_\epsilon}}\,\right|\,\right)\,\bigg|\,d\sigma(x)\\
 & = 
\displaystyle\int_{\partial\Omega} \left|\,[\,g'(\gamma(u_{\epsilon}+ \xi_x)) - g'(\gamma(u_{\epsilon}))\,]\gamma(\,w_1(x)-w_2(x)\,)\,\gamma(\Phi)\gamma\left(\,\left|\,\displaystyle\frac{J_{\partial\Omega} {h_\epsilon}}{J {h_\epsilon}}\,\right|\,\right)\, \right|\,d\sigma(x) \\ 
 & \leq K_1 \,\left\{\displaystyle\int_{\partial\Omega} [\,(g'(\gamma(u_{\epsilon}+ \xi_x)) - g'(\gamma(u_{\epsilon})))]^{\,p}[\gamma(w_1(x)-w_2(x))]^{p}\left[\gamma\left(\left|\displaystyle\frac{J_{\partial\Omega} {h_\epsilon}}{J {h_\epsilon}}\right|\right)
\right]^{p}d\sigma(x)\right\}^\frac{1}{p}||\Phi|| {_{X^\frac{1}{2}}}\,, \\
 & \leq K_1  K_ 2 \, \mu_{\epsilon} \,  \,\left\{\displaystyle\int_{\partial\Omega} [\,(g'(\gamma(u_{\epsilon}+ \xi_x)) - g'(\gamma(u_{\epsilon})))]^{\,p}d\sigma(x)\right\}^\frac{1}{p} 
 \|w_1 - w_2\|_{ X^{\eta}   }||\Phi|| {_{X^\frac{1}{2}}}\,, \\
\end{align*}
where 
  $\mu_{\epsilon} = \left|\displaystyle\frac{J_{\partial\Omega} {h_\epsilon}}{J {h_\epsilon}}\right|$ is  bounded, uniformly in  $\epsilon$
  and  $ w_1(x) \leq \xi_x \leq w_2(x)$ or $ w_2(x) \leq \xi_x \leq w_1(x)$.

Now the integrand above is bounded by   $2^{\,{p}}||\,g'\,||_\infty^{\,{p}}||$  and goes to $0$ a.e. as $\rho \to 0$, since $||\,w_1\,||_{X^{\eta}}\leq \rho$,  $||\,w_2\,||_{X^{\eta}}\leq \rho$  and  $w_1(x)\leq \xi_x\leq w_2(x)$. 
 Thus, the integral goes to $0$ by Lebesgue's dominated convergence Theorem.
 \eproof
 
 We are now in a position to prove the main result of this section

\begin{teo}
Assume the hypotheses of Theorem \ref{equicon} hold.  Then the family of attractors
     \{$\mathcal{A}_{\epsilon }\,{|}\, 0 \leq \epsilon \leq \epsilon_0\}$, of
     the problem (\ref{abstract_scale}), whose existence is guaranteed by
     Theorem \ref{exatr}  is lower semicontinuous in   $X^\eta$.
  \end{teo}
\proof

The system generated by
(\ref{abstract_scale}) is gradient for any $\epsilon$ and its equilibria are all hyperbolic for $\epsilon $ in a neighborhood of $0$. Also, 
 the equilibria are continuous in $\epsilon$ by Theorem \ref{equicon}, the  linearisation is  continuous
 in $\epsilon$ as shown during  the proof of  Theorem \ref{equicon} and  the local unstable manifolds of the equilibria are continuous
 in $\epsilon$, by Theorem \ref{manifcont}. The result follows then from \cite{PP}, Theorem 3.10 .  
\eproof 


\end{document}